\documentclass{article} 

\usepackage{authblk} 
\usepackage[a4paper,top=3cm,bottom=3cm,left=3cm,right=3cm]{geometry} 
\usepackage{mathtools}
\usepackage[english]{babel}

\usepackage[T1]{fontenc}
\usepackage{amssymb}
\usepackage{amsmath}
\usepackage{amsthm}      
\usepackage{tikz-cd}
\usepackage{xfrac}
\usepackage{faktor}
\usepackage{mathrsfs}
\usepackage{graphicx}
\usepackage{mathabx}
\usepackage{float}
\usepackage{xcolor}         
\usepackage{listings}          
\usepackage{hyperref}     
\usepackage[normalem]{ulem}

\usepackage{bbm}
\usepackage{enumitem}
\usepackage{comment}

\usepackage{amssymb}
\usepackage[refpage]{nomencl}
\makenomenclature

\usepackage{etoolbox}
\renewcommand\nomgroup[1]{%
  \item[\bfseries
  \ifstrequal{#1}{M}{Maps}{%
  \ifstrequal{#1}{S}{Sets}{%
  \ifstrequal{#1}{R}{Representations}{
  \ifstrequal{#1}{O}{Other symbols}{}}}}%
]}

\providetoggle{nomsort}
\settoggle{nomsort}{true} 

\makeatletter
\iftoggle{nomsort}{%
    \let\old@@@nomenclature=\@@@nomenclature        
        \newcounter{@nomcount} \setcounter{@nomcount}{0}%
        \renewcommand\the@nomcount{\two@digits{\value{@nomcount}}}
        \def\@@@nomenclature[#1]#2#3{
          \addtocounter{@nomcount}{1}%
        \def\@tempa{#2}\def\@tempb{#3}%
          \protected@write\@nomenclaturefile{}%
          {\string\nomenclatureentry{\the@nomcount\nom@verb\@tempa @[{\nom@verb\@tempa}]%
          \begingroup\nom@verb\@tempb\protect\nomeqref{\theequation}%
          |nompageref}{\thepage}}%
          \endgroup
          \@esphack}%
      }{}
\makeatother

\newtheorem{theorem}{Theorem}[section]
\newtheorem{defi}[theorem]{Definition}
\newtheorem{cor}[theorem]{Corollary}
\newtheorem{lemma}[theorem]{Lemma}
\newtheorem{prop}[theorem]{Proposition}

\theoremstyle{definition}
\newtheorem{rmk}[theorem]{Remark}

\makeatletter
\newsavebox{\@brx}
\DeclareMathOperator{\cind}{c-ind}
\DeclareMathOperator{\pind}{i}

\newcommand{\llangle}[1][]{\savebox{\@brx}{\(\m@th{#1\langle}\)}%
  \mathopen{\copy\@brx\mkern2mu\kern-0.9\wd\@brx\usebox{\@brx}}}
\newcommand{\rrangle}[1][]{\savebox{\@brx}{\(\m@th{#1\rangle}\)}%
  \mathclose{\copy\@brx\mkern2mu\kern-0.9\wd\@brx\usebox{\@brx}}}
\makeatother

\newcommand{\of}[1]{\mathcal{O}_{#1}}
\newcommand{\pf}[1]{\mathfrak{p}_{#1}}

\author{Elena Collacciani \\ \href{mailto:elena.collacciani@phd.unipd.it}{elena.collacciani@phd.unipd.it} }

\affil{Padova University}

\title{A Reduction over finite fields of the tame local Langlands correspondence for $SL_N$ }

\begin{document}

\maketitle

\begin{abstract}
  \noindent We establish the conjecture formulated by Vogan in \cite{Voganslide} for \( SL_n \) . Specifically, we construct a surjection from the set of irreducible representations of \( SL_n(k) \), where \( k \) is a finite field, to the inertia equivalence classes of tame Langlands parameters for \( SL_n(F) \), where \( F \) is a \( p \)-adic field with residue field \( k \). Additionally, we provide a parametrization of the fibers of this surjection and examine its compatibility with the Local Langlands Correspondence for \( SL_n \). 
  
  \noindent This work extends several results previously established for \( GL_n \) in \cite{Macdonald, Zink2} to the context of \( SL_n \).

\end{abstract}


\section*{Introduction}
The primary goal of this paper is to prove, in the $SL_n$ case, the conjecture formulated by Vogan in \cite{Voganslide}. This conjecture aims to establish a Langlands parameterization for finite groups of Lie type.
\\
\\
Let $F$ be a non-archimedean local field of characteristic $0$ with residue field $k_F$. Let $\mathbb{G}$ be a connected, reductive, linear algebraic group defined over $F$ and $F$-split. We denote the group of $F$-points of $\mathbb{G}$ by $\mathbf{G} = \mathbb{G}(F)$.

A tame Langlands parameter for $\mathbf{G}$ is as a pair $(\rho, E)$, where $\rho$ is a continuous homomorphism from the Weil group of $F$ into $\mathbb{G}^{\vee}(\mathbb{C})$, the complex dual group of $\mathbb{G}$ \cite[Section 2.1]{Bor}, and $E$ is a nilpotent element in the complex Lie algebra of $\mathbb{G}^{\vee}$ \cite[Section 8.1]{Bor}, satisfying some compatibility relations.
The (tame) Langlands correspondence is, conjecturally, a surjective map $\mathbf{\mathcal{L}}$ with finite fibers, called \textit{L-packets}, from the set of isomorphism classes of (depth-$0$) smooth admissible representations of $\mathbf{G}$ onto the set of equivalence classes of (tame) Langlands parameters for $\mathbf{G}$. Within an $L$-packet, the representations are conjecturally parametrized by irreducible representations of a - slightly modified - component group associated with the centralizer of the tame Langlands parameter. A detailed exposition of these conjectures lies beyond the scope of this paper. We refer the interested readers to \cite{Kal}
\\
\\
The finite group of Lie type $\overline{\mathbf{G}} = \mathbb{G}(k_F)$ given by the group of the $k_F$-points of $\mathbb{G}$ can be seen as a quotient of the maximal hyperspecial subgroup $\mathbf{K} = \mathbb{G}(\of{F})$. Vogan's conjecture proposes a procedure to derive a suitable class of parameters for $\overline{\mathbf{G}}$ from the tame Langlands parameters for $\mathbf{G}$. Furthermore, for each parameter, the corresponding "packet" should be parametrized by the irreducible representations of a component group derived from the parameter, akin to the \( p \)-adic case.

The most comprehensive and general parameterization for irreducible representations of finite groups of Lie type is due to Lusztig \cite{Orange,Lus88}. While it resembles a Langlands parameterization in some respects, it differs in some important aspects. According to Lusztig's parameterization, the representations of \(\overline{\mathbf{G}}\) are classified by:
\begin{enumerate}[label=L.\arabic*]
    \item A semisimple conjugacy class $s$ in $\mathbb{G}^{\vee}(k_F)$, the dual group of $\mathbb{G}$ over $k_F$;
    \item A special nilpotent orbit \(E\) in the Lie algebra of the centralizer of \(s\);
    \item An element of the set \(M(A):=\faktor{\{ (x,\psi) \mid x\in A, \, \psi \text{ irreducible representation of }C_A(x)\}}{\sim_{A}},\) 
    where \(\sim_A\) denotes \(A\)-conjugation, and \(A\) is the so-called canonical quotient - a specific quotient of the component group of  of the centralizer of $(s,E)$ in $\mathbb{G}^{\vee}(\overline{k_F})$, the dual group of $\mathbb{G}$ over the algebraic closure of $k_F$.
\end{enumerate}

Vogan's conjecture builds on the expectation that Lusztig's parameterization can be modified to resemble other Langlands parameterizations. Specifically, the desired parameters should involve:
\begin{enumerate}[label=V.\arabic*]
    \item A semisimple conjugacy class \(s\) in \(\mathbb{G}^{\vee}(\mathbb{C})\), the dual group of $\mathbb{G}$ over the complex field  \label{1};
    \item A nilpotent orbit \(E\) in the Lie algebra of the centralizer of \(s\) (relaxing the "special" condition)\label{2};
    \item An irreducible representation of the component group of a centralizer in $\mathbb{G}^{\vee}(\mathbb{C})$ depending on \(s\) and \(E\)\label{3}. 
\end{enumerate}

\hspace{2mm}

\noindent The construction in \cite{Macdonald} for \(GL_N\) serves as an inspiring model. In this case, the representations of \(GL_N(k_F)\) are parametrized by "$I_F$-equivalence classes" of tame Langlands parameters for $GL_N(F)$. An $I_F$-equivalence class is, roughly speaking, defined by the restriction of \((\rho, E)\) - a tame Langlands parameter for \(GL_N(F)\) - to the inertia subgroup. This corresponds to the desiderata above: \ref{1} is given by a generator of the image through $\rho$ of the inertia subgroup , and \ref{2} is given by \(E\). In the case of \(GL_N\), the connectedness of centralizers renders \ref{3} irrelevant.

Vogan's approach generalizes Macdonald's construction to other groups by refining the definition of $I_F$-equivalence classes. The key idea is to preserve partial information about $\rho(Fr)$, the image of a Frobenius element of $F$ under the map $\rho$. Specifically, he proposes to retain in the parameterization the coset of the connected component of the stabilizer in $\mathbb{G}^{\vee}(\mathbb{C})$ of $\rho(I_F)$, the image of the inertia subgroup under $\rho$, determined by $\rho(Fr)$.

The conjecture states that there exists a surjection from the irreducible representations of $\overline{\mathbf{G}}$ to the set of these equivalence classes, with the representations within each fiber of this map parametrized by irreducible representations of the component group of the stabilizer in $\mathbb{G}^{\vee}(\mathbb{C})$ of the equivalence class relative to the fiber.

Additionally, Vogan’s conjecture suggests a compatibility between the proposed parameterization of irreducible representations of $\overline{\mathbf{G}}$ and the Langlands parameterization of depth-$0$ representations of $\mathbf{G}$. Specifically, if $\pi$ is an  irreducible representation of $\overline{\mathbf{G}}$ associated to the $I_F$-equivalence class of the tame Langlands parameter $(\rho,E)$ then the inflation of $\pi$ to the maximal compact subgroup $\mathbf{K}$ of $\mathbf{G}$ should appear as a $\mathbf{K}$-subrepresentation of some representation in the \(L\)-packet corresponding to \((\rho, E)\).
Furthermore, there exists a natural map $\iota$ between the component group of the centralizer in $\mathbb{G}^{\vee}(\mathbb{C})$ of $(\rho,E)$, whose irreducible representations parametrize the $L$-packet for $(\rho, E)$, and the component group of the centralizer in $\mathbb{G}^{\vee}(\mathbb{C})$ of the $I_F$-equivalence class of $(\rho,E)$, whose irreducible representations parametrize the fiber over the $I_F$-equivalence class of $(\rho, E)$.
If $\pi$ is parametrized by an irreducible representation $\zeta$ of the latter component group, the representation in the $L$-packet corresponding to $(\rho,E)$ containing the inflation of $\pi$ as $\mathbf{K}$-subrepresentation should be precisely the one parametrized by the irreducible representation $\zeta\circ\iota$.

This statement is already established for \(GL_N\) by \cite[Appendix A]{Zink2}, where Silberger and Zink prove the compatibility of Langlands parameterization with Macdonald’s parameterization. In this case, the connectedness of centralizers ensures that both \(L\)-packets and Macdonald’s packets contain only a single representation.
\\
\\
To the best of the author’s knowledge, Macdonald’s work \cite{Macdonald} is the only example where a parameterization of the type envisioned by Vogan has been explicitly studied. The objective of this paper is to demonstrate that Vogan’s conjecture also holds for \(SL_N\).

The first section contains a review of key results for \(GL_N\) from \cite{Macdonald} and \cite{Zink2}, which provide a foundation for the subsequent analysis. In the second section, we focus on the \(SL_N\) case,  beginning with a concise review of \cite{Gel}, where the existence of a local Langlands correspondence \(\mathcal{L}'\) and the parameterization of \(L\)-packets by component groups are established.

We start addressing Vogan's conjecture in Section \ref{Voganconj}. In \eqref{MV} we define the "Macdonald–Vogan correspondence" \(\mathcal{M}'_N\) , which yields a parametrization of the irreducible representations of $SL_N(k_F)$ as envisioned by Vogan. Theorem \ref{fiber2} establishes the parameterization of its fibers by irreducible representations of the expected component groups. Although this parameterization is not canonical, it becomes so upon fixing a "base point" representation. The proof of  Theorem \ref{fiber2} draws inspiration from \cite[Section 4]{Gel}, adapting arguments to the finite field case.

The remainder of the paper is dedicated to proving the compatibility of parameterizations provided by $\mathcal{L}_N'$ and $\mathcal{M}_N'$ respectively. The key results are Theorem \ref{compatibilitysln}, establishing compatibility of \(L\)-packets for \(SL_N(F)\) with the fibers of the Macdonald-Vogan correspondence for \(SL_N(k_F)\), and Theorem \ref{finalcomp}, confirming Vogan’s predictions regarding the interplay between  the parametrization by irreducible representations of component groups and \({K'_N}\)-restrictions.

Our results rely heavily on the strong connection between the representation theories of $SL_N$ and the one of $GL_N$, the latter of  which benefits from an especially extensive body of results. However, we aim to extend this study to other split reductive groups in future work.

\section*{Notation}
We summarize here some general notation we will use in the course of the paper. 
\begin{itemize}
\item Let $F$ \nomenclature{$F$}{Non-archimedean local field} be a non-archimedean local field of characteristic $0$ with residue field $k_F$ \nomenclature{$k_F$}{Residue field of $F$}, a finite field with $q$ elements. Let $\of{F}$ \nomenclature{$\of{F}$}{Ring of integer of $F$} be the valuation ring in $F$ and $\pf{F}$ \nomenclature{$\pf{F}$}{Maximal ideal of $\of{F}$} be the maximal ideal of $\of{F}$. we denote by $\overline{k_F}$ the algebraic closure of $k_F$

We denote by $W_{F}$  \nomenclature{$W_F$}{Weil group of $F$} the Weil group of $F$. We denote by $I_F$ \nomenclature{$I_F$}{Inertia subgroup of $W_F$} and $P_F$ \nomenclature{$P_F$}{Wild ramification subgroup of $W_F$} respectively the inertia and the wild ramification subgroups of $W_F$, and by $Fr$ \nomenclature{$Fr$}{Frobenius element of $W_F$} a Frobenius element of $W_F$.

\item Let $H$ be a group. We denote by $\widehat{H}$ the character group of $H$, sometimes with $H^{\wedge}$ for notation reasons. With a slight abuse of notation, when $H=F^*$ or $H=W_F$, we write $\widehat{F^*}$ \nomenclature{ $\widehat{F^*}$}{Group of tame characters of $F$} and $\widehat{W_F}$ \nomenclature{ $\widehat{W_F}$}{Group of tame characters of $F$} for the groups of tamely ramified characters of $F^*$ and $W_F$ respectively.

If $H$ acts on a set $X$, for any $x\in X$ we write  $stab_H(x)$ for the subgroup of $H$ fixing $x$
\item Let $H$ be a group and let $K\leq H$. If $\pi$ is a representation of $H$, we write $Res^{H}_{K}$ \nomenclature{$Res$}{Restriction}, or sometimes $\pi|_{K}$, for the restriction of $\pi$ to $K$.\\
If $\gamma$ is a representation of $K$, we write $ind_K^H\gamma$ \nomenclature{ $ind$}{Induction} for the induced representation and $\cind_{K}^{H}\gamma$  \nomenclature{ $\cind$}{Compact induction} for the compactly induced representation.  For any $h\in H$, we write ${}^h\gamma$ for the representation of ${}^hK=hKh^{-1}$ defined by ${}^h\gamma(x)=\gamma(h^{-1}xh)$ for any $x\in{}^hK$.

Let $N$ be a normal subgroup of $H$. If $\pi$ is a representation of $H$, we denote by $\pi^{N}$ \nomenclature{ $(-)^{N}$}{$N$-fixed points} the $H/N$ representation on the $N$-fixed subspace of $\pi$. If $\gamma$ is a representation of $H/N$, we denote by $Infl_{\faktor{H}{N}}^{H}\gamma$  \nomenclature{$Infl$}{Inflation} the inflated representation.  

Let $\pi$ be a representation $H$ and $\pi'$ is a subrepresentation of $\pi$, not necessarily irreducible. We say that $\pi'$ has multiplicity $1$ in $\pi$ if $Hom_{H}(\pi',\pi)=Hom_{H}(\pi',\pi')$.

\item If $\mathbb{H}$ is a reductive linear algebraic group defined over a field $\mathbb{F}$, we denote by $\mathbb{H}(\mathbb{F})$ the $\mathbb{F}$-points of $\mathbb{H}$.
\item Assume $\mathbb{H}$ to be defined over $\mathbb{C}$. If $X\subseteq \mathbb{H}(\mathbb{C})$, we write $C_{\mathbb{H}(\mathbb{C})}(X)=\{h\in\mathbb{H}(\mathbb{C}) | hXh^{-1}=X\}$  \nomenclature{$C_{G}(*)$}{Centralizer in $G$ of $*$} for the stabilizer of $X$ in $\mathbb{H}$ with respect to the conjugation action. If $f$ is a map to $\mathbb{H}(\mathbb{C})$, we denote $C_{\mathbb{H}(\mathbb{C})}(f)=\{h\in\mathbb{H}(\mathbb{C}) | hxh^{-1}=x \ \text{ for any } \ x\in Im(f)\}$ the centralizer of the image of $f$. If $x\in\mathfrak{h}$, the complex Lie algebra of $\mathbb{H}(\mathbb{C})$, we set $C_{\mathbb{H}(\mathbb{C})}(x)=\{h\in\mathbb{H}(\mathbb{C}) | Ad(h)x=x\}$. We denote $C_{\mathbb{H}(\mathbb{C})}(f,X,x):= C_{\mathbb{H}(\mathbb{C})}(f)\cap C_{\mathbb{H}(\mathbb{C})}(X)\cap C_{\mathbb{H}(\mathbb{C})}(x)$. Similarly, if $\mathcal{Y}$ is a collection of subsets of $\mathbb{H}(\mathbb{C})$, maps to $\mathbb{H}(\mathbb{C})$ and elements of $\mathfrak{h}$, we set $C_{\mathbb{H}(\mathbb{C})}(\mathcal{Y}):=\bigcap_{\textit{y}\in\mathcal{Y}}C_{\mathbb{H}(\mathbb{C})}(\textit{y})$ and $C^0_{\mathbb{H}}(\mathcal{Y})$ for the connected component of $C_{\mathbb{H}}(\mathcal{Y})$. We denote by $A_{\mathbb{H}}(\mathcal{Y}):=\faktor{C_{\mathbb{H}}(\mathcal{Y})}{C_{\mathbb{H}}(\mathcal{Y})^0}$ \nomenclature{$A_{G}(*)$}{Component group of the centralizer in $G$ of $*$} the component group of $C_{\mathbb{H}}(\mathcal{Y})$. 

\item Let $\mathbb{G}$ be a reductive, connected linear algebraic group defined over $F$ and $F$-split. 

Let $\mathbf{G}=\mathbb{G}(F)$ be the group of $F$-points of $\mathbb{G}$. Then $\mathbf{K}$ denotes the hyperspecial maximal compact subgroup $=\mathbb{G}(\of{F})$ of  $\mathbf{G}$ and we denote by
$\mathbf{K}^+$ its pro-unipotent radical, setting $\overline{\mathbf{G}}:=\faktor{\mathbf{K}}{\mathbf{K}^+}\cong\mathbb{G}(k_F)$.

Let $\mathbf{\Omega}_0$ be the set of isomorphism classes of  irreducible admissible representations of $\mathbf{G}$ of depth $0$, and let $\overline{\mathbf{\Omega}}$ be the set of the isomorphism classes of irreducible  representations of $\overline{\mathbf{G}}$.
We write the parahoric restriction functor for $\mathbf{G}$ as
\begin{align}\label{parres} 
\mathcal{P}^{\mathbf{G}}_{\overline{\mathbf{G}}}:\mathbf{\Omega}_0&\rightarrow \overline{\mathbf{\Omega}}\\
    \pi&\mapsto (Res^{\mathbf{G}}_{\mathbf{K}}\pi)^{\mathbb{K^+}}\nonumber   
    \end{align}
\nomenclature{$\mathcal{P}$}{Parahoric restriction}
Let $\mathbb{G}^{\vee}$ denote the dual group of $\mathbb{G}$ and $\mathfrak{g}^{\vee}$ the complex Lie algebra of $\mathbb{G}^{\vee}$. A tame Langlands parameter for $\mathbf{G}$ is a pair $(\rho, E)$, where $\rho:W_F\rightarrow \mathbb{G}^{\vee}(\mathbb{C})$ is a continuous homomorphism such that $\rho(Fr)$ is semisimple and $\rho|_{P_F}=1$, and $E\in\mathfrak{g}^{\vee} $ is a nilpotent element satisfying $Ad(\rho(w))(E)=\|w\| E$ for any $w\in W_F$, where $\|w\|$ is the norm of $W_F$. The group $\mathbb{G}^{\vee}(\mathbb{C})$ acts on the set of Langlands parameters by the adjoint action, and two parmeters are equivalent if they are in the same $\mathbb{G}^{\vee}(\mathbb{C})$-orbit. We denote by $\mathbf{\Phi}_0$ the set of  of tame Langlands parameters for $\mathbf{G}$. 

In the course of the paper, we will be concerned with $\mathbb{G}=GL_n$ and $\mathbb{G}=SL_n$, with $n\in\mathbb{N}_{\geq1}$. When we specialize to $\mathbb{G}=GL_n$, we write $G_n:=\mathbf{G}=GL_n(F)$, and in general we denote all the objects defined in the previous paragraph specialized to the $GL_n$-case with the same symbols used before but in normal (not-bold) font and with a subscript $n$ (e.g. $K_n$ in place of $\mathbf{K}$,  or $(\Omega_n)_0$ in place of $\mathbf{\Omega}_0$).\\
The complex dual group of $GL_N$ is $GL_N(\mathbb{C})$, so a tame Langlands parameter for $G_N$ is a pair $(\rho, E)$ where $\rho:W_F\mapsto GL_N(\mathbb{C})$ is a continuous group morphism trivial on $P_F$ and with semisimple image, and $E\in\mathfrak{gl}_N(\mathbb{C})$ is a nilpotent element satisfying $Ad(\rho(w))(E)=\|w\| E$ for any $w\in W_F$. The set $(\Phi_N)_0$ is the set of equivalence classes of  tame Langlands parameters for $G_N$, where the equivalence is given by $GL_N(\mathbb{C})$ adjoint action.

Similarly, when we specialize to $\mathbb{G}=SL_n$, we write $G'_n:=\mathbf{G}=SL_n(F)$, and in general we denote all the objects defined in the previous paragraph specialized to this case with the same symbols used before but in normal (not-bold) font, with a subscript $n$ and a superscript $'$ (e.g. $K'_n$ in place of $\mathbf{K}$,  or $(\Omega'_n)_0$ in place of $\mathbf{\Omega}_0$). \\
The complex dual group of $SL_N$ is $PGL_N(\mathbb{C})$, so a tame Langlands parameter for $G'_N$ is a pair $(\rho, E)$ where $\rho:W_F\mapsto PGL_N(\mathbb{C})$  is a continuous group morphism trivial on $P_F$ and with semisimple image, and $E\in\mathfrak{sl}_N(\mathbb{C})$ is a nilpotent element satisfying $Ad(\rho(w))(E)=\|w\| E$ for any $w\in W_F$.  The set $(\Phi'_N)_0$ is the set of equivalence classes of  tame Langlands parameters for $G'_N$, where the equivalence is given by $PGL_N(\mathbb{C})$ adjoint action.

\nomenclature{$G_n$}{The group $GL_n(F)$}
\nomenclature{$K_n$}{The subgroup $GL_n(\of{F})\leq G_n$}
\nomenclature{$K^+_n$}{The nilpotent radical of $K_n$}
\nomenclature{$\overline{G}_n$}{The group $\faktor{K_n}{K^+_n}\cong GL_n(k_F)$}
\nomenclature{$(\Omega_n)_0$}{Smooth admissible irreducible depth-$0$ representations of $G_n$}
\nomenclature{$(\overline{\Omega_n})_0$}{Irreducible representations of $\overline{G}_n$}
\nomenclature{$(\Phi_n)_0$}{Tame Langlands parameters for $G_n$}
\nomenclature{$G'_n$}{The group $SL_n(F)$}
\nomenclature{$K'_n$}{The subgroup $SL_n(\of{F})\leq G'_n$}
\nomenclature{${K'}^+_n$}{The nilpotent radical of $K'_n$}
\nomenclature{$\overline{G'}_n$}{The group $\faktor{K'_n}{{K'}^+_n}\cong SL_n(k_F)$}
\nomenclature{$(\Omega'_n)_0$}{Smooth admissible irreducible depth-$0$ representations of $G'_n$}
\nomenclature{$(\overline{\Omega'_n})_0$}{Irreducible representations of $\overline{G'}_n$}
\nomenclature{$(\Phi'_n)_0$}{Tame Langlands parameters for $G'_n$}

\item We denote the normalized parabolic induction functor by $\pind$ \nomenclature{$\pind$}{Parabolic induction}. We will need it just for the general linear group, so we use a notation tailored on the situations in which we will use it.  Let $n\in\mathbb{N}$ and $(m_1,\dots, m_k)$ be a composition of $n$, and let $\pi_j$ (respectively $\overline{\pi}_j$) be a representation of $G_{m_j}$ (respectively $\overline{G}_{m_j}$) for $j\in\{1\dots k\}$. Then we write $\pind_{\prod_{j=1}^{k}G_{m_j}}^{G_n} \bigboxtimes_{j=1}^{k}\pi_j$ (respectively $\pind_{\prod_{j=1}^{k}\overline{G}_{m_j}}^{G_n} \bigboxtimes_{j=1}^{k}\overline{\pi}_j$) for the parabolic induction to $G_n$ (respectively to $\overline{G}_n$) from the standard Levi subgroup sitting diagonally in $G_n$ (respectively $\overline{G}_n$) with blocks of the prescribed dimensions, where the inflation is performed through the parabolic subgroup of block upper triangular matrices containing the aforementioned Levi subgroup. 
\end{itemize}
For the readers' convenience, an index of the symbols used has been added at the end of the paper.
\newpage
\section{Compatibility of Langlands and Macdonald Correspondences for $GL_N$}
In this section, we reformulate the result in \cite[Appendix A]{Zink2}, that was stated in terms of the "reduction of the tempered type" map, \cite{Zink}. We will be interested in the depth-$0$ case, where the map can be described just in terms of the parahoric restriction functor \eqref{parres}.In the level $0$ case, the latter coincides  with the functors denoted by "$k_{max}$" and "$T_{k,\lambda}$" in \cite{Zink}. We adopt this easier description because  the "reduction of the tempered type" map is not yet available for a general split reductive group. The statements in terms of parahoric restriction will make the extension of the results to the $SL_N$ case more straightforward.

\subsection{The Head of Parahoric restriction}\label{Gln}
Let $n\in\mathbb{N}$.  Let  $(C_n)_0\subseteq (\Omega_n)_0$ be the set of isomorphism classes of depth-0 supercuspidal irreducible representations of $G_n$, and let ${C}_0=\bigsqcup ({C}_n)_0$. For any $\tau\in {C}_0$, we say that $\tau$ has degree $n$ if $\tau\in({C}_n)_0$, and we write $d(\tau)=n$. We say that $\tau_1, \tau_2\in (C_n)_0$  are equivalent up to unramified twist if there exists an unramified character $\chi$ of  $G_n$ such that $\tau_2=\tau_1\otimes\chi$. Let  $(\mathcal{C}_n)_0$ be the set of equivalence classes of depth-0 supercuspidal representations of ${G_n}$ up to unramified twists, and let $\mathcal{C}_0=\bigsqcup (\mathcal{C}_n)_0$.

Let $\overline{\mathcal{C}_n}$ be the set of isomorphism classes of cuspidal representation of $\overline{G}_n$ and $\overline{\mathcal{C}}=\bigsqcup \overline{\mathcal{C}_n}$. For any $\overline{\tau}\in \overline{\mathcal{C}}$, we say that $\overline{\tau}$ has degree $n$ if $\overline{\tau}\in(\mathcal{C}_n)_0$, and we write $d(\overline{\tau})=n$.\\

For any $[\tau]\in (\mathcal{C}_n)_0$  there exists a $\overline{\tau}\in\overline{\mathcal{C}_n}$ such that for any $\sigma\in[\tau]$ with central character $\omega_\sigma\in\widehat{F^*}$ it holds $\sigma=\cind_{F^*K_n}^{G_n} (\omega_\sigma \otimes \overline{\tau})$. By \cite[Lemma 5.5]{Zink}, 
\begin{equation}
\mathcal{P}^{G_n}_{\overline{G}_n}\sigma=\overline{\tau}.
\end{equation}
Moreover any $\sigma\in (C_n)_0$ obtained by compact induction from $\overline{\tau}\in\overline{\mathcal{C}_n}$ belongs to the same  class in $(\mathcal{C}_n)_0$. Hence the parahoric restriction functor yields a degree-preserving bijection
\begin{align}\label{parrescusp}
    \mathcal{C}_0&\rightarrow \overline{\mathcal{C}}\\
    [\tau] &\mapsto\overline{\tau}:=\mathcal{P}^{G_N}_{\overline{G}_N}\tau,\nonumber
    \end{align}
See also \cite[Fact 1, Appendix A 1.2]{Zink2}.

Let $Par(n)$ be the set of the partitions of $n$, and $Par=\bigsqcup (Par(n))$. For $\lambda\in Par(n)$ we write $|\lambda|=n$.

A partition valued function on $\mathcal{C}_0$ is a finite support function $\Lambda:\mathcal{C}_0\mapsto Par$. Its degree is given by
 \[d(\Lambda):=\sum_{[\tau]\in \mathcal{C}_0}d(\tau)|\Lambda([\tau])|.\]
The dominance order $\leq_{dom}$ on partitions induces a partial order on partition valued functions on $\mathcal{C}_0$: we say that $\Lambda_1\geq\Lambda_2$ if and only if  for any $[\tau]\in\mathcal{C}_0$ it holds $|\Lambda_1([\tau])|=|\Lambda_2([\tau])|$ and $\Lambda_1([\tau])\leq_{dom}\Lambda_2([\tau])$.

 To any representation $\pi \in (\Omega_N)_0$ is associated a partition valued function on $\mathcal{C}_0$ of degree $N$ as follows. Following \cite{Zel}, a segment $\Delta=\Delta(\tau,l)$ with $l\in\mathbb{N}$ and $\tau\in C_0$  is a subset of $C_0$ of the form
\[\Delta(\tau,l):=\{\tau\otimes |det(\cdot)|^i, \ i=0\dots l-1\}\]
and by Zelevisky's classification (\cite[Theorem 6.1, Section 9.1]{Zel}) to any  $\pi \in \Omega_0$ is assigned a set of segments $\{\Delta(\tau_1,l_1), \dots, \Delta(\tau_k,l_k)\}$. For any $[\tau]\in\mathcal{C}_0$, we define $J([\tau]):=\{i\in \{1,\dots, k\} \ | \ \tau_i\in [\tau] \}$. Then the partition valued function associated to $\pi$ is given by
\begin{equation}\label{deflambdapi}
\Lambda_{\pi}([\tau]):=(l_i | i\in J([\tau])),
\end{equation}
\nomenclature{$\Lambda_{\pi}$}{Partition valued function associated to $\pi$ \eqref{deflambdapi}}
that is $\Lambda_{\pi}([\tau])$ the partition given by the lengths of the segments in the Zelevinky's classification of $\pi$ containing supercuspidal representations in the unramified twists class $[\tau]$.  \\

A partition valued function on $\overline{\mathcal{C}}$ is a function $M:\overline{\mathcal{C}}\mapsto Par$. Its degree is given by
 \[d(M):=\sum_{\overline\tau\in\overline{ \mathcal{C}}}d(\overline\tau)|M(\overline{\tau})|.\] 

There is a bijection between the set of partition valued functions on $\overline{\mathcal{C}}$ of degree $N$ and $\overline{\Omega_N}$,  obtained as follows: for any $\overline{\tau}\in \overline{\mathcal{C}_n}$ and for any $m\in\mathbb{N}$, the endomorphism algebra $End_{\overline{G}_{nm}}(\pind_{\overline{G}_{n}^{m}}^{\overline{G}_{nm}}(\bigboxtimes_{i=1}^m\overline{\tau}))$ is a finite Hecke algebra relative to the symmetric group of $m$ elements $\mathbb{S}_{m}$, and therefore it is isomorphic to the group algebra $\mathbb{C}[\mathbb{S}_{m}]$ \cite{HLHecke}. So the irreducible constituents of $\pind_{\overline{G}_{n}^{m}}^{\overline{G}_{nm}}(\bigboxtimes_{i=1}^m\overline{\tau})$ 
are in bijection with simple (right) modules over $End_{\overline{G}_{nm}}(\pind_{\overline{G}_{n}^{m}}^{\overline{G}_{nm}}(\bigboxtimes_{i=1}^m\overline{\tau}))$, that in turn are in one to one correspondence to the irreducible representations of $\mathbb{S}_{m}$, which are parametrized by partitions of $m$. We normalize the parametrization in such a way that the partition $(m)$ corresponds to the sign representation. For any $\lambda\in Par(m)$ we write $\overline{\tau}^{\lambda}$ 
\nomenclature{$\overline{\tau}^{\lambda}$}{Irreducible representation with cuspidal support $\bigotimes_{i=1}^m{\overline\tau}$ corresponding to $\lambda\vdash m$ } 
for the irreducible constituent of $\pind_{\overline{G}_{n}^{m}}^{\overline{G}_{nm}}(\bigboxtimes_{i=1}^m\overline{\tau})$ such that $Hom_{\overline{G}_{nm}}(\pind_{\overline{G}_{n}^{m}}^{\overline{G}_{nm}}(\bigboxtimes_{i=1}^m\overline{\tau}), \overline{\tau}^{\lambda})$ as $End_{\overline{G}_{nm}}(\pind_{\overline{G}_{n}^{m}}^{\overline{G}_{nm}}(\bigboxtimes_{i=1}^m\overline{\tau}))$-right module is the simple module corresponding to $\lambda$.
The bijection between partition valued functions on $\overline{\mathcal{C}}$ of degree $N$ and $\overline{\Omega_N}$ is given by
\begin{align}\label{sigmabar}
\overline{\sigma}:\{ M:\overline{\mathcal{C}}\mapsto Par \ | \ d(M)=N\}&\longrightarrow \overline{\Omega_N}\\
M &\mapsto \overline{\sigma}_{M}:=\pind_{\prod_{\overline{\tau}\in\overline{\mathcal{C}}}\overline{G}_{d(\overline{\tau})|M(\overline{\tau})|}}^{\overline{G}_N}(\bigboxtimes_{\overline{\tau}\in\overline{\mathcal{C}}}\overline{\tau}^{M(\overline{\tau})}) \nonumber
\end{align}
where $ \overline{\sigma}_{M}$ is irreducible by \cite{HLcomp}, see also \cite[pag. 3]{Macdonald}.
\nomenclature{$\overline{\sigma}_{\overline{\Lambda}}$}{Irreducible representation of $\overline{G}_n$ corresponding to $\overline{\Lambda}$ in \eqref{sigmabar}}

The parahoric restriction functor induces  through the bijection \eqref{parrescusp} the bijection 
\begin{align}\label{deflambdabar}
\{ \Lambda:\mathcal{C}_0\mapsto Par \ | \ d(\Lambda)=N\}&\rightarrow \{ M:\overline{\mathcal{C}}\mapsto Par \ | \ d(M)=N\}\\
\Lambda &\mapsto \overline{\Lambda}:(\overline{\tau}\mapsto\Lambda([\tau)])\nonumber
\end{align} 

The equations \eqref{deflambdapi}, \eqref{deflambdabar} and \eqref{sigmabar} together yield a map
\begin{align}\label{head1}
(\Omega_{N})_0&\longrightarrow \overline{\Omega_N}\\
\pi &\mapsto \overline{\sigma}_{\overline{\Lambda_{\pi}}}\nonumber
\end{align} 

\nomenclature{$\overline{\Lambda}$}{Partition valued function on $\overline{\mathcal{C}}$ associated to $\overline{\Lambda}$ by \eqref{deflambdabar}}
The following theorem is a restatement, in the depth-$0$ case, of \cite[Proposition 6.2]{Zink}, see also \cite[A.1.3, fact 4]{Zink2}. 

\begin{theorem} \label{prephead}\cite[Proposition 6.2]{Zink}
    For any $\pi\in \Omega_N$, the representation $\overline{\sigma}_{\overline{\Lambda_{\pi}}}$ of $\overline{G_N}$ obtained via \eqref{head1} is an irreducible constituent of multiplicity $1$ of $\mathcal{P}^{G_N}_{\overline{G}_N}\pi$. If $\overline{\sigma}_{\overline{\Lambda}}$ for $\Lambda$ a partition valued function of degree $N$ is an irreducible constituent of $\mathcal{P}^{G_N}_{\overline{G}_N}\pi$, then ${\Lambda}\leq{\Lambda_{\pi}}$. If $\pi$ is tempered, then $\sigma_{\overline{\Lambda}}$ is an irreducible constituent of $\mathcal{P}^{G_N}_{\overline{G}_N}\pi$ for any partition valued function $\Lambda$ of degree $N$ with ${\Lambda}\leq {\Lambda_{\pi}}$.
\end{theorem}

\begin{rmk}
In \cite[Proposition 6.2]{Zink} the partition valued functions are denoted by $\mathcal{P}$ rather than $\Lambda$, and the representations denoted there by $\sigma_{\mathcal{P}}(\lambda)$ are the inflations to $K_N$ of the representations of $\overline{G}_n$ that we denote by $\overline{\sigma}_{\overline{\Lambda}}$ in Theorem \ref{prephead}. Moreover the relation between the partition valued functions $\Lambda_{\pi}$ and the representation $\pi$ is expressed in the reference by $\pi\in im(Q_{\mathcal{\pi}})$.

We chose our notation to be more similar to the one used in \cite[A.1.3, fact 4]{Zink2}, where the authors deal with the depth $0$ case directly. The representations denoted there by $\sigma_{\Lambda}$ are the inflations to $K_N$ of the representations of $\overline{G}_n$ that we denote by $\overline{\sigma}_{\overline{\Lambda}}$ in Theorem \ref{prephead}.
\end{rmk}

\begin{defi}
In the notation  above, the Head of the parahoric restriction of $\pi$ is the irreducible representation of $\overline{G}_N$ given by
\begin{equation}\label{Head}
\mathcal{HP}^{G_N}_{\overline{G}_N}\pi:= \overline{\sigma}_{\overline{\Lambda_{\pi}}}=\pind_{\prod_{\overline{\sigma}\in\overline{\mathcal{C}}}\overline{G}_{d(\overline{\sigma})|\overline{\Lambda_{\pi}}(\overline{\sigma})|}}^{\overline{G}_N}(\bigboxtimes_{\overline{\sigma}\in\overline{\mathcal{C}}}\overline{\sigma}^{\overline{\Lambda_{\pi}}(\overline{\sigma})}.)
\end{equation}
where $\overline{\sigma}_{\overline{\Lambda_{\pi}}}$ is the image of $\pi$ through the map \eqref{head1}.
\nomenclature{$\mathcal{HP}^{G_n}_{\overline{G}_n}$}{Head of parahoric restriction for $G_n$ \eqref{Head}}
\end{defi}
By Theorem \ref{prephead}, the representation $\mathcal{HP}^{G_N}_{\overline{G}_N}\pi$ is an irreducible constituent of $\mathcal{P}^{G_N}_{\overline{G}_N}\pi$ of multiplicity $1$.\\

\subsection{Compatibility of Langlands and Macdonald correspondence}
\paragraph{Level-$0$ Local Langlands correspondence}
The local Langlands correspondence for $GL$ is a family of bijections between Langlands parameters for $G_n$ and irreducible representations of $G_n$ for any $n\in\mathbb{N}$ satisfying a list of properties. The existence (and uniqueness)  of the  local Langlands correspondence for $GL$ was established independently in \cite{HT, Hen, Sc}. Each bijection in the local Langlands correspondence restricts to a bijection between tame Langlands parameters for $G_n$ and depth-$0$ irreducible representations of $G_n$. We denote the bijections by
\begin{align}\label{Langland}
&\mathcal{L}_n: (\Omega_n)_0\rightarrow (\Phi_n)_0&\text{for $n\in\mathbb{N}$.}
 \end{align}
 \nomenclature{$\mathcal{L}_n$}{Tame local Langlands correspondence for $G_n$ \eqref{Langland}}
  In the level $0$ case, the bijections $\mathcal{L}_n$ are explicitly described in \cite{BH1, BH2, BH3}. 
 
 By \cite[Proposition 3.2]{BH1}, the map $\mathcal{L}_1$ is the restriction to tame characters of the bijection of local class field theory induced by the Artin reciprocity map \cite[Section 29]{BKGL2}, and for any  tame character $\chi\in \widehat{F^*}$ and any $\pi\in(\Omega_n)_0$ it holds
 \begin{equation}\label{detequiL}
     \mathcal{L}_N(\chi\circ det \otimes \pi) = \mathcal{L}_1(\chi)\mathcal{L}_N(\pi)
 \end{equation}
where the product on the right-hand side is the multiplication of a tame character of $W_F$ with the homomorphism $\rho:W_F\mapsto GL_N(\mathbb{C})$.
\paragraph{Macdonald correspondence}
 Following \cite{Macdonald}, two tame Langlands parameters for $G_n$, say  $(\rho_i, E_i)$, for $i=1,2$ are called $I_F$-equivalent if there exists $A\in GL_n(\mathbb{C})$ such that
\begin{align}
A(\rho_1|_{I_F}) A^{-1}&=\rho_2|_{I_F}, \label{aGln}\tag{a}  \\
 Ad(A)E_1&=E_2. \label{bGln}\tag{b}
\end{align}
We denote by $\faktor{(\Phi_n)_0}{\sim_{I_F}}$ the set of $I_F$-equivalence classes in $(\Phi_n)_0$, and we write $(\rho, E)_{I_F}\in\faktor{(\Phi_N)_0}{\sim_{I_F}}$ for the $I_F$-equivalence class of $(\rho, E)\in(\Phi_N)_0$. The Macdonald correspondence is the family of bijections  
\begin{align}\label{macdonaldgen}
    &\mathcal{M}_n: \overline{\Omega_n} \rightarrow \faktor{({\Phi_n})_0}{\sim_{I_F}}&\text{for any $n\in\mathbb{N}$}
\end{align}
\nomenclature{$\mathcal{M}_n$}{Macdonald correspondence for $G_n$ \eqref{macdonaldgen}}
 constructed in \cite{Macdonald}. See also \cite{Aub} for a recent description and characterization of \eqref{macdonaldgen}.\\
 By \cite[Proposition 1.2]{Macdonald}, the map $\mathcal{M}_1$ is the bijection induced by local class field theory between Frobenius stable characters in $\widehat{\faktor{I_F}{P_F}}$ and $\widehat{k_F^*}$. For any Frobenius stable character $\chi\in\widehat{\faktor{I_F}{P_F}}$ and any $\pi\in \overline{\Omega_n}$ it holds
  \begin{equation}\label{detequiM}
     \mathcal{M}_n((\chi\circ det )\otimes \pi) = \mathcal{M}_1(\chi)\mathcal{M}_n(\pi)
 \end{equation}
by \cite[Proposition 1.3]{Macdonald}, where the product on the right hand side is induced from the multiplication action of $\widehat{W_F}$ on $\Phi_n$.

\begin{rmk}\label{rmkIFGln}
 If $(\rho_i, E_i)\in (\Phi_N)_0$ for $i=1,2$ are such that $(\rho_1, E_1)\sim_{I_F}(\rho_2, E_2)$, by definition there exists $A\in GL_n(\mathbb{C})$ satisfying \eqref{aGln},\eqref{bGln}.Then
\begin{equation}\label{cGLn}
A( \rho_1(Fr) C^0_{GL_n(\mathbb{C})}(\rho_1|_{I_F}))A^{-1}=\rho_2(Fr) C^0_{GL_n(\mathbb{C})}(\rho_2|_{I_F}).\tag{c}
\end{equation}
Centralizers in $GL_n(\mathbb{C})$ are connected,  so $C^0_{GL_n(\mathbb{C})}(\rho_i|_{I_F})=C_{GL_n(\mathbb{C})}(\rho_i|_{I_F})$. Moreover by condition \eqref{aGln}
\[A C_{GL_n(\mathbb{C})}(\rho_1|_{I_F})A^{-1}=C_{GL_n(\mathbb{C})}(A(\rho_1|_{I_F})A^{-1})=C_{GL_n(\mathbb{C})}(\rho_2|_{I_F}).\]
 Therefore \eqref{cGLn} holds if and only if 
\[\rho_2(Fr)^{-1}A \rho_1(Fr)A^{-1}\in C^0_{GL_n(\mathbb{C})}(\rho_2|_{I_F}).\]
The $\rho_i$ are tame representations, so they are representations of $\faktor{W_F}{P_F}=\faktor{I_F}{P_F}\rtimes \langle Fr\rangle$, where $Fr$ acts by powering to the $q$ every element in $\faktor{I_F}{P_F}$. Therefore using condition \eqref{aGln} twice
\begin{align*}
\rho_2(Fr)^{-1}A \rho_1(Fr)A^{-1}\rho_2|_{I_F} A \rho_1(Fr)^{-1}A^{-1} \rho_2(Fr)=\\
\rho_2(Fr)^{-1}A \rho_1(Fr) \rho_1|_{I_F}\rho_1(Fr)^{-1}A^{-1} \rho_2(Fr)=\\
\rho_2(Fr)^{-1}A  (\rho_1|_{I_F})^q A^{-1} \rho_2(Fr)=\\
\rho_2(Fr)^{-1}  (A\rho_1|_{I_F}A^{-1})^q \rho_2(Fr)=\\
\rho_2(Fr)^{-1}(\rho_2|_{I_F})^q \rho_2(Fr)= \rho_2|_{I_F}.
\end{align*}
Hence $\rho_2(Fr)^{-1}A \rho_1(Fr)A^{-1}\in C^0_{GL_n(\mathbb{C})}(\rho_2|_{I_F})$, so $I_F$-equivalent representations in $(\Phi_n)_0$ always satisfy \eqref{cGLn}
\end{rmk}

\paragraph{Compatibility of Langlands and Macdonald correspondence} \begin{theorem}\label{compatibility} \cite[A.2, Proposition 2]{Zink2}
For any $N\in\mathbb{N}$ the following diagram commutes:
\begin{equation}
\begin{tikzcd}
(\Omega_N)_0 \arrow[d, "\mathcal{HP}^{G_N}_{\overline{G}_N}"] \arrow[r, "\mathcal{L}_N"] &(\Phi_N)_0  \arrow[d, "\faktor{}{\sim_{I_F}}"] \\
 \overline{\Omega_N}\arrow[r, "\mathcal{M}_N"]  &\faktor{(\Phi_N)_0}{\sim_{I_F}}                                  
\end{tikzcd}
\end{equation}
Where $\mathcal{HP}^{G_N}_{\overline{G}_N}$ is the Head of the parahoric restriction map defined in \eqref{Head}
\end{theorem}

\section{Local Langlands correspondence for $SL_N(F)$}
We recall the construction of the local Langlands correspondence for $SL_N(F)$ as carried out in \cite{Gel}. The two "Working Hypothesis" assumed in \cite{Gel} have been confirmed:  the existence of a local Langlands correspondence for $G_N$ was established independently in \cite{HT, Hen, Sc}, and the restrictions of irreducible representations of $G_N$  to $G'_N$ have been proved to be multiplicity free in \cite{Tadic}.\\
The group $\widehat{F^*}$ acts on $(\Omega_N)_0$, with a character $\chi$ acting by $\widetilde{\pi}\mapsto\widetilde{\pi} \otimes(\chi\circ det)$ for $\widetilde{\pi}\in(\Omega_N)_0$. The group $\widehat{W_F}$ acts on $(\Phi_N)_0$ by multiplication. Identifying $\widehat{F^*}$ and  $\widehat{W_F}$  by the bijection of local class field theory, the local Langlands correspondence for $G_N$ is equivariant with respect to these group actions by \eqref{detequiL}.
Therefore, denoting the orbit-sets of the action of the tame character groups by $\faktor{(\Phi_N)_0}{\widehat{W_F}}$ and $\faktor{(\Omega_N)_0}{\widehat{F^*}}$, the local Langlands correspondence induces a bijection
\[\overline{\mathcal{L}}_N: \faktor{(\Omega_N)_0}{\widehat{F^*}}\ \rightarrow\faktor{(\Phi_N)_0}{\widehat{W_F}}.\]


Let $\eta:GL_N(\mathbb{C})\rightarrow PGL_N(\mathbb{C})$ \nomenclature{$\eta$}{Projection $GL_N(\mathbb{C})\rightarrow PGL_N(\mathbb{C})$ }be the natural projection. Then the map
\begin{align}
    \eta^{*}:\faktor{(\Phi_N)_0}{\widehat{W_F}}\rightarrow (\Phi'_N)_0\\
    \widehat{W_F}(\rho,E)\mapsto (\eta\circ \rho , E)\nonumber
\end{align}
is a bijection.\\
For any $\widetilde{\pi}\in (\Omega_N)_0$, the representation $Res^{G_N}_{G'_N}\widetilde{\pi}$ is a direct sum of finitely many mutually inequivalent representations of $G'_N$. Moreover, for  $\widetilde{\pi}_1, \widetilde{\pi}_2\in(\Omega_N)_0$ their restrictions to $G'_N$ are either equal, and that happens if and only if $\widetilde{\pi}_1=\widetilde{\pi}_2\otimes \chi\circ det$ for some  $\chi\in \widehat{F^*}$, or they do not have any irreducible constituent in common. Let $(\Omega'_N)_0$ be the set of isomorphism classes of irreducible smooth admissible representations of $G'_N$ with supercuspidal support of depth $0$. For any $\pi\in (\Omega'_N)_0$, let $R(\pi)\in {(\Omega_N)_0}$ \nomenclature{$R(\pi)$}{Representative of $\mathcal{R}(\pi)$} be an irreducible representation of $G_N$ containing $\pi$ as $G'_N$ subrepresentation. Then the map 
\begin{align}\label{R}
\mathcal{R}: (\Omega'_N)_0&\rightarrow \faktor{(\Omega_N)_0}{\widehat{F^*}}\\
\pi&\mapsto \widehat{F^*}R(\pi)\nonumber
\end{align}
\nomenclature{$\mathcal{R}$}{The map \eqref{R}}
does not depend on the choice of $R(\pi)$ and it is a surjection. We will usually denote by $R(\pi)$ a representative of the orbit $\mathcal{R}(\pi)$

The local Langlands correspondence for $G'_N$ is the surjection given by the composition
\begin{equation}\label{L'}
    \mathcal{L}'_N:(\Omega'_N)_0\xrightarrow{\mathcal{R}} \faktor{(\Omega_N)_0}{\widehat{F^*}}\xrightarrow {\overline{\mathcal{L}}_N}\faktor{(\Phi_N)_0}{\widehat{W_F}}\xrightarrow{\eta^{*}} (\Phi'_N)_0.
\end{equation}
\nomenclature{$\mathcal{L}'_n$}{The tame local Langlands correspondence for $G'_n$ \eqref{L'}}
The fibers of this surjection are called $L$-packets. For any $(\rho, E)\in (\Phi'_N)_0$,  the  character group of the component group $A_{PGL_N(\mathbb{C})}(\rho, E)$ has a canonical simply transitive action on the $L$-packet ${ \mathcal{L}'_N}^{-1}(\rho,E)$ \cite[Theorem 4.3]{Gel}.

\section{The Macdonald-Vogan correspondence for $SL_N(k_F)$}\label{Voganconj}

\subsection{The $I_F$-equivalence classes}
\begin{defi}\label{IFequivsln}
The elements $(\rho_2, E_2),(\rho_1, E_1)\in (\Phi'_N)_0$  are $I_F$-equivalent if there exists $A\in PGL_N(\mathbb{C})$ such that
\begin{align}
    A\rho_1|_{I_F}A^{-1}&=\rho_2 \label{a}\tag{a'}\\
   Ad(A) E_1&= E_2\label{b}\tag{b'}\\
 A \big(\rho_1(Fr) C^0_{PGL_N(\mathbb{C})}(\rho_1|_{I_F})\big)A^{-1}&=\rho_2(Fr) C^0_{PGL_N(\mathbb{C})}(\rho_2|_{I_F}).\label{c}\tag{c'}
\end{align}
We denote by $\sim_{I_F}$ the $I_F$-equivalence relation. We write $(\rho, E)_{I_F}\in\faktor{(\Phi'_N)_0}{\sim_{I_F}}$ for the $I_F$-equivalence class of $(\rho, E)\in(\Phi'_N)_0$  
\nomenclature{$\sim_{I_F}$}{$I_F$-equivalence relation}
\end{defi}
\begin{rmk}\label{rmkIFGln2}
In view of  Remark \eqref{rmkIFGln}, this definition of $I_F$-equivalence is actually analogous to the one given for $G_N$.
\end{rmk}

The group $\widehat{k_F^*}$ acts on $\faktor{(\Phi_N)_0}{\sim_{I_F}}$ as follows: the group $\widehat{k_F^*}$ can be identified with the group of Frobenius stable characters of $\faktor{I_F}{P_F}$ via class field theory. Any Frobenius stable character $\chi$ of $\faktor{I_F}{P_F}$ can be extended trivially to $W_F$, and we still denote by $\chi$ the extended character. Then $\chi$ acts on $(\Phi_N)_0$ by $(\rho, E)\mapsto(\rho\otimes\chi, E)$. Since $(\rho_1, E_1)\sim_{I_F}(\rho_2, E_2)$ implies $(\rho_1\otimes\chi, E_1)\sim_{I_F}(\rho_2\otimes\chi, E_2)$, this defines an action of $\widehat{k_F^*}$ on $\faktor{(\Phi_N)_0}{\sim_{I_F}}$.

\begin{lemma}
Let $\eta:GL_N(\mathbb{C})\rightarrow PGL_N(\mathbb{C})$ be the natural projection. The assignment $(\rho, E)\mapsto (\eta\circ \rho,E) $, for $(\rho, E)\in (\Phi_N)_0$, induces a well defined bijection
 \begin{align}\label{vogbij}
   \eta^{*}:\faktor{\big(\faktor{(\Phi_N)_0}{\sim_{I_F}}\big)}{\widehat{k_F^*}}\rightarrow \faktor{(\Phi'_N)_0}{\sim_{I_F}},
 \end{align}
where the action of $\widehat{k_F^*}$ on $\faktor{(\Phi_N)_0}{\sim_{I_F}}$ is the one defined above.
 \begin{proof}
     We first show that the map $(\rho, E)\mapsto (\eta\circ \rho,E) $  induces a well defined map
     \begin{align}\label{IFeqcomp}
  \faktor{(\Phi_N)_0}{\sim_{I_F}}&\rightarrow \faktor{(\Phi'_N)_0}{\sim_{I_F}}.
 \end{align}
 
     If  $(\rho_i, E_i)\in (\Phi_N)_0$, for $i=1,2$, are such that $(\rho_1, E_1)\sim_{I_F}(\rho_2, E_2)$, by Remark \ref{rmkIFGln} there exists $A\in GL_N(\mathbb{C})$ satisfying conditions \eqref{aGln}. \eqref{bGln}, \eqref{cGLn}.     
     Therefore applying $\eta$ gives 
     \begin{align*}
         \eta(A)(\eta\circ\rho_1|_{I_F})\eta(A^{-1})&=\eta\circ A\rho_1|_{I_F}A^{-1}=\eta\circ \rho_2|_{I_F} &\text{condition \eqref{a}}\\
         Ad(\eta(A)) E_1&= E_2 &\text{condition \eqref{b}}.
     \end{align*}
    We check condition \eqref{c}. Observe that $\eta (C_{GL_N(\mathbb{C})}(\rho_i|_{I_F}))=C^0_{PGL_N(\mathbb{C})}(\eta\circ \rho_i|_{I_F})$ for $i=1,2$. Then there holds
     \begin{align*}
     \eta(A)(\eta\circ\rho_1(Fr)C^0_{PGL_N(\mathbb{C})}(\eta\circ \rho_1|_{I_F}))\eta(A^{-1})=\\
     \eta(A)\eta (\rho_1(Fr)C_{GL_N(\mathbb{C})}(\rho_1|_{I_F}))\eta(A^{-1})=\\
     \eta(A\rho_1(Fr)C_{GL_N(\mathbb{C})}(\rho_1|_{I_F})A^{-1})=\\
     \eta(
     \rho_2(Fr) C_{GL_N(\mathbb{C})}(\rho_2|_{I_F}))=\\
     (\eta\circ\rho_2(Fr)) C^0_{PGL_N(\mathbb{C})}(\eta\circ \rho_2|_{I_F}).
     \end{align*}
So $(\eta\circ \rho_1, E_1)\sim_{I_F}(\eta\circ\rho_2, E_2)$ in $(\Phi'_N)_0$.

The map \eqref{IFeqcomp} is constant on  $\widehat{k_F^*}$-orbits because $(\eta\circ (\rho\otimes\chi),E)= (\eta\circ \rho,E)$ for any $(\rho, E)\in (\Phi_N)_0$, hence $\eta^*$ is well defined.\\

We show that it is injective. Let $(\rho_i, E_i)\in (\Phi_N)_0$ for $i=1,2$ be such that $(\eta\circ\rho_1, E_1)\sim_{I_F} (\eta\circ\rho_2, E_2)$ in $(\Phi'_N)_0$. Then there exists $A\in PGL_N(\mathbb{C})$ satisfying conditions \eqref{a},\eqref{b} and \eqref{c}

Let $\tilde{A}\in GL_N(\mathbb{C})$ be such that $\eta(\tilde{A})=A$. Then $Ad(\tilde{A})E_1=E_2$, and $\tilde{A}\rho_1|_{I_F}\tilde{A}^{-1}=\chi\rho_2|_{I_F}$ for some  character $\chi$ of $I_F$. Since $\rho_1$ and $\rho_2$ are tame, $\chi$ is a character of $\faktor{I_{F}}{P_F}$, so  $(\rho_1,E_1)\sim_{I_F}(\rho_2\otimes\chi, E_2)$ in $(\Phi_N)_0$. In order to prove injectivity, it is therefore sufficient to show that $\chi$ is Frobenius stable. For $i=1,2$ it holds
\[\rho_i(Fr)\rho_i|_{I_F}\rho_i(Fr^{-1})=(\rho_i|_{I_F})^q.\]
It follows that
\begin{equation}\label{chirat1}
    \tilde{A}\rho_1(Fr)(\rho_1|_{I_F})\rho_1(Fr^{-1})\tilde{A}^{-1}= \tilde{A}(\rho_1|_{I_F})^q\tilde{A}^{-1}=(\tilde{A}(\rho_1|_{I_F})\tilde{A}^{-1})^q=(\chi\rho_2|_{I_F})^q=\chi^q(\rho_2|_{I_F})^q.
\end{equation}
Moreover lifting to $GL_N(\mathbb{C})$ the relation
\[A( \eta\circ \rho_1(Fr) C^0_{PGL_N(\mathbb{C})}(\eta\circ \rho_1|_{I_F}))A^{-1}=\eta\circ \rho_2(Fr) C^0_{PGL_N(\mathbb{C})}(\eta\circ\rho_2|_{I_F})\]
we obtain 
\[\tilde{A} \rho_1(Fr) C_{GL_N(\mathbb{C})}(\rho_1|_{I_F})\tilde{A}^{-1}=\rho_2(Fr) C_{GL_N(\mathbb{C})}(\rho_2|_{I_F}).\]
In particular there exists a $c\in C_{GL_N(\mathbb{C})}(\rho_2|_{I_F})$ such that 
\[\tilde{A} \rho_1(Fr)\tilde{A}^{-1}=\rho_2(Fr) c\]
Therefore
\begin{align}\label{chirat2}
    \tilde{A}\rho_1(Fr)(\rho_1|_{I_F})\rho_1(Fr^{-1})\tilde{A}^{-1}= 
    \tilde{A}\rho_1(Fr) \tilde{A} \tilde{A}^{-1}(\rho_1|_{I_F} )\tilde{A}^{-1} \tilde{A}\rho_1(Fr^{-1})\tilde{A}^{-1}= \nonumber \\
    \rho_2(Fr) c \chi (\rho_{2}|_{I_F} )c^{-1}\rho_2(Fr^{-1})=\chi \rho_2(Fr) (\rho_{2}|_{I_F}) \rho_2(Fr^{-1}) =
    \chi(\rho_{2}|_{I_F})^q    
\end{align}
Comparing equations \eqref{chirat1} and \eqref{chirat2} we get
\[\chi^q (\rho_2|_{I_F})^q=  \tilde{A}\rho_1(Fr)(\rho_1|_{I_F})\rho_1(Fr^{-1})\tilde{A}^{-1}=  \chi(\rho_{2}|_{I_F})^q,\]
giving $\chi^q=\chi$. Hence, $\chi$ is a Frobenius stable character of $\faktor{I_F}{P_F}$, i.e. it is an element of $\widehat{k_F^*}$. In other words, the $I_F$-equivalence classes in $(\Phi_N)_0$ of $(\rho_1,E_1)$ and $(\rho_2, E_2)$ are in the same $ \widehat{k_F^*}$-orbit. \\

It remains to show that \eqref{vogbij} is surjective. Let $(\overline{\rho},E)\in(\Phi'_N)_0$. The existence of  $(\rho,E)\in(\Phi_N)_0$ such that $(\overline{\rho},E)=(\eta\circ\rho,E)$ amounts to the existence of a semisimple Weil representation $\rho:W_F\rightarrow GL_N(\mathbb{C})$ that is a lift of $\overline{\rho}$, which was proved in \cite{Hen_repWF}.
\end{proof}
\end{lemma}
\subsection{The Macdonald-Vogan correspondence}
The group $\widehat{k_F^*}$ acts on $\overline{\Omega_N}$ and on $\faktor{(\Phi_N)_0}{\sim_{I_F}}$,  with a character $\chi$ acting respectively by $\pi\mapsto\pi\otimes(\chi\circ det)$ for $\pi\in \overline{\Omega_N}$, and by $(\rho, E)_{I_F}\mapsto(\rho\otimes\chi, E)_{I_F}$ for $(\rho,E)\in\faktor{(\Phi_N)_0}{\sim_{I_F}}$. By \eqref{detequiM}, the Macdonald correspondence \eqref{macdonaldgen}
is equivariant with respect to these actions of $\widehat{k_F^*}$, so it induces a bijection between orbit sets
\[\overline{\mathcal{M}}_N:  \faktor{\overline{\Omega_N}}{\widehat{k_F^*}}\rightarrow\faktor{\big(\faktor{(\Phi_N)_0}{\sim_{I_F}}\big)}{\widehat{k_F^*}}.\]

Since $\overline{G'}_N$ is normal in $\overline{G}_N$, the latter acts on $\overline{\Omega'_N}$ by conjugation, inducing an action of $\faktor{\overline{G}_N}{\overline{G'}_N}\cong k_F^*$  on $\overline{\Omega'_N}$. 

For any $\pi\in \overline{\Omega'_N}$ let $\overline{R}(\pi)\in \overline{\Omega_N}$  be such that $\pi$ is a subrepresentation of $Res^{\overline{G}_N}_{\overline{G'}_N}\overline{R}(\pi)$. \nomenclature{$\overline{R}(\pi)$}{Representative of $\overline{\mathcal{R}}(\pi)$} 
From Clifford theory it follows that the map
\begin{align}\label{bijR}
\faktor{\overline{\Omega'_N}}{k^*_F}&\rightarrow\faktor{\overline{\Omega_N}}{\widehat{k^*_F}}\\
k^*_F\pi&\mapsto \widehat{k^*_F}\overline{R}(\pi) \nonumber
\end{align}
is a bijection \cite[Proposition 5.1]{Macdonald}.  
In particular, the map
\begin{align}\label{Rbar}
\overline{\mathcal{R}}: \overline{\Omega'_N}&\rightarrow\faktor{\overline{\Omega_N}}{\widehat{k^*_F}}\\
\pi&\mapsto \widehat{k^*_F}\overline{R}(\pi) \nonumber
\end{align}
\nomenclature{$\overline{\mathcal{R}}$}{The map \eqref{Rbar}}
is a well-defined surjection. In the following, we will denote by $\overline{R}(\pi)$ a representative of the orbit $\overline{\mathcal{R}}(\pi)$.

\begin{defi}  
We call the Macdonald-Vogan correspondence the surjection
    \begin{equation}\label{MV}
\mathcal{M}'_N:
\overline{\Omega'_N}\xrightarrow{\overline{\mathcal{R}}} \faktor{\overline{\Omega_N}}{\widehat{k^*_F}}\xrightarrow {\overline{\mathcal{M}}_N}\faktor{\big(\faktor{(\Phi_N)_0}{\sim_{I_F}}\big)}{\widehat{k_F^*}}\xrightarrow{\eta^{*}} \faktor{(\Phi_N)_0}{\sim_{I_F}}
\end{equation}
\nomenclature{$\mathcal{M}'_n$}{The Macdonald-Vogan correspondence \eqref{MV}}
\end{defi}

\section{Parametrization of the fibers of $\mathcal{M}'_N$}
The goal of this Section is Theorem \ref{fiber2}, which yields a parameterization of each fiber of $\mathcal{M}'_N$ in terms of irreducible representations of a component group. 

\begin{lemma}\label{fiber0}
For $\pi\in\overline{\Omega'_N}$ let $\overline{R}(\pi)\in \overline{\Omega_N}$ be such that $\pi$ is an irreducible $G'_N$ subrepresentation of $\overline{R}(\pi)$. 
There is a canonical isomorphism
 \[\faktor{k_F^*}{Stab_{k_F^*}(\pi)}=(Stab_{\widehat{k^*_F}}(\overline{R}(\pi))) ^{\wedge}.\]
 where $(Stab_{\widehat{k^*_F}}(\overline{R}(\pi)) )^{\wedge}$ denotes the character group of $Stab_{\widehat{k^*_F}}(\overline{R}(\pi))$.
\begin{proof}
We show that $Stab_{\widehat{k^*_F}}
(\overline{R}(\pi))=\big(\faktor{k_F^*}{Stab_{k_F^*}(\pi)}\big)^{\wedge}$, where the equality makes sense viewing $\chi\in\big(\faktor{k_F^*}{Stab_{k_F^*}(\pi)}\big)^{\wedge}$ as a character in $\widehat{k_F^*}$ that is trivial on $Stab_{k_F^*}(\pi)$. The statement will then follow by duality.   

We start by proving that $\big(\faktor{k_F^*}{Stab_{k_F^*}(\pi)}\big)^{\wedge}\leq Stab_{\widehat{k^*_F}}(\overline{R}(\pi))$.

Let $ \chi\in\widehat{k_F^*}$ such that $\chi|_{Stab_{k_F^*}(\pi)}=1$. 
The determinant map gives an isomorphism $\faktor{\overline{G}_N}{\overline{G'}_N}\cong  k^*_F$ and the action of $k^*_F$ on $\overline{\Omega_N'}$ is induced by the conjugation action of $\overline{G}_N$. It follows that the determinant map  induces a group isomorphism $\faktor{\overline{G}_N}{Stab_{\overline{G}_N}(\pi)}\cong\faktor{k_F^*}{Stab_{k_F^*}(\pi)}$,
 that in turn induces the group isomorphism $\chi\mapsto\chi\circ \det$ between characters of $k_F^*$ that are trivial on $Stab_{k_F^*}(\pi)$ and characters of $\overline{G}_N$ that are trivial on $Stab_{\overline{G}_N}(\pi)$.
 
 Let $\{g_1\dots g_k\}$ be a set of representatives of the cosets of $Stab_{G_N}(\pi)$ in $\overline{G}_N$. Then $Res^{\overline{G}_N}_{\overline{G'}_N}\overline{R}(\pi)=\bigoplus_{i=1}^k {}^{g_i}\pi$. Let $A:=\bigoplus_{i=1}^k \chi(\det(g_i^{-1}))$, that is, the linear map acting as the scalar $\chi(\det(g_i^{-1}))$ on the subspace on which $Res^{\overline{G}_N}_{\overline{G'}_N}\overline{R}(\pi)$ acts as $^{g_i}\pi$. Then 
\[A\overline{R}(\pi)A^{-1}=\overline{R}(\pi)\otimes (\chi\circ\det).\] 
Hence
\[\overline{R}(\pi)\otimes \chi\circ \det \cong \overline{R}(\pi)\]
that is, $\chi\in Stab_{\widehat{k_F^*}}(\overline{R}(\pi))$.  So $\big(\faktor{k_F^*}{Stab_{k_F^*}(\pi)}\big)^{\wedge}\leq Stab_{\widehat{k^*_F}}(\overline{R}(\pi))$.

By Clifford theory, \cite[Proposition 5.1]{Macdonald}
\begin{equation*}
    |k^*_F \pi|=|Stab_{\widehat{k^*_F}}(\overline{R}(\pi))|.
\end{equation*}
 Therefore
 \[|\big(\faktor{k_F^*}{Stab_{k_F^*}(\pi)}\big)^{\wedge}|=|\faktor{k_F^*}{Stab_{k_F^*}(\pi)}|=|k^*_F \pi|=|Stab_{\widehat{k^*_F}}(\overline{R}(\pi))|\]
 and hence, by cardinality,
 \begin{equation}\label{isomorphismfiber1}
\big(\faktor{k_F^*}{Stab_{k_F^*}(\pi)}\big)^{\wedge}=Stab_{\widehat{k^*_F}}(\overline{R}(\pi)). 
 \end{equation}

 Since $\faktor{k_F^*}{Stab_{k_F^*}(\pi)}$ is a finite abelian group, it is canonically isomorphic to its double character group, hence dualizing \eqref{isomorphismfiber1} we get
 \[\faktor{k_F^*}{Stab_{k_F^*}(\pi)}=(Stab_{\widehat{k^*_F}}(\overline{R}(\pi))) ^{\wedge}.\]
\end{proof}
\end{lemma}

\begin{lemma}\label{fiber1}
Let $(\rho,E)_{I_F}\in  \faktor{(\Phi'_N)_0}{\sim_{I_F}}$ and $(\tilde{\rho}, E)_{I_F}\in \faktor{(\Phi_N)_0}{\sim_{I_F}}$ be such that $(\eta\circ\tilde{\rho}, E)_{I_F}=(\rho, E)_{I_F}$, and let
\[Stab_{\widehat{k^*_F}}((\tilde{\rho}, E)_{I_F})=\{\chi\in\widehat{k^*_F} \ | \ (\chi\otimes \tilde{\rho},E)_{I_F}=(\tilde{\rho}, E)_{I_F}\}.\]
 Then the character group of ${Stab_{\widehat{k^*_F}}((\tilde{\rho}, E))}$ acts simply transitively on ${\mathcal{M}'_N}^{-1}((\rho,E)_{I_F})$.
\begin{proof}
Set $\tilde{\pi}_{(\tilde{\rho},E)}:=\mathcal{M}_N^{-1}((\tilde{\rho},E)_{I_F})$.
By construction $\mathcal{M}'_N=(\eta^{*})\circ \overline{\mathcal{M}}_N\circ \overline{\mathcal{R}}$.  The maps $\overline{\mathcal{M}}_N$ and $(\eta^{*})$ are bijections,  so $\widehat{k^*_F} \tilde{\pi}_{(\tilde{\rho}, E)}$ is the unique $\widehat{k^*_F}$-orbit in $\overline{\Omega_N}$ satisfying $(\rho,E)=(\eta^{*})\circ\overline{\mathcal{M}}_N(\widehat{k^*_F} \tilde{\pi}_{(\tilde{\rho}, E)})$. Then 
 \[{\mathcal{M}'_N}^{-1}((\rho,E)_{I_F})= \overline{\mathcal{R}}^{-1}(\widehat{k^*_F} \tilde{\pi}_{(\tilde{\rho}, E)}).\]
 The map $\overline{\mathcal{R}}$ factorizes as
\begin{align*}
  \overline{\Omega'_N}\xrightarrow{/k^*_F} &\faktor{\overline{\Omega'_N}}{k^*_F}\rightarrow\faktor{\overline{\Omega_N}}{\widehat{k^*_F}}
\end{align*}
where the last map is the bijection defined in \eqref{bijR}. Therefore if $\pi_{(\rho,E)}\in \overline{\Omega'_N}$ satisfies $\overline{R}(\pi_{(\rho,E)})=\widehat{k^*_F} \tilde{\pi}_{(\tilde{\rho},E)}$, then
\[\overline{\mathcal{R}}^{-1}(\widehat{k^*_F} \tilde{\pi}_{(\tilde{\rho},E)})=k^*_F \pi_{(\rho,E)}\]
It follows that the group $\faktor{k_F^*}{Stab_{k_F^*}(\pi_{(\rho, E)})}$ has a canonical simply transitive action on ${\mathcal{M}'_N}^{-1}((\rho, E)_{I_F})$. 
Therefore by Lemma \eqref{fiber0}, the character group $(Stab_{\widehat{k^*_F}}(\tilde{\pi}_{(\tilde{\rho},E)}))^{\wedge}$ has a simply transitive action on ${\mathcal{M}'_N}^{-1}((\rho, E)_{I_F})$.
 
Since the map $\mathcal{M}_N$ is compatible with the $\widehat{k^*_F}$-actions, we have
\[Stab_{\widehat{k^*_F}}(\tilde{\pi}_{(\tilde{\rho},E)})= Stab_{\widehat{k^*_F}}((\tilde{\rho},E)_{I_F})=\{\chi\in\widehat{k^*_F} \ | \ (\chi\otimes \tilde{\rho},E)_{I_F}= (\tilde{\rho},E)_{I_F}\}\]
and dualizing

\[(Stab_{\widehat{k^*_F}}(\tilde{\pi}_{(\tilde{\rho},E)}))^{\wedge}= (Stab_{\widehat{k^*_F}}((\tilde{\rho},E)_{I_F}))^{\wedge}=\{\chi\in\widehat{k^*_F} \ | \ (\chi\otimes \tilde{\rho},E)_{I_F}=(\tilde{\rho},E)_{I_F}\}^{\wedge}\]
\end{proof}
\end{lemma}

\begin{theorem}\label{fiber2}
Let  $(\rho,E)_{I_F}\in  \faktor{(\Phi'_N)_0}{\sim_{I_F}}$. 
The character group $A_{PGL_N(\mathbb{C})}(\rho|_{I_F},\rho(Fr)C^0_{PGL_N(\mathbb{C})}(\rho|_{I_F}),E)^{\wedge}$ is a finite abelian group and it acts simply transitively on  
${\mathcal{M}'_N}^{-1}((\rho,E)_{I_F})$

\begin{proof}
Let $(\tilde{\rho},E)_{I_F}\in \faktor{(\Phi_N)_0}{\sim_{I_F}}$ be such that $(\eta\circ\tilde{\rho},E)_{I_F}=(\rho,E)_{I_F}$. By Lemma \ref{fiber1} it is enough to provide a group isomorphism between
$A_{PGL_N(\mathbb{C})}(\rho|_{I_F},\rho(Fr)C^0_{PGL_N(\mathbb{C})}(\rho|_{I_F}),E)$ and $Stab_{\widehat{k^*_F}}((\tilde{\rho},E)_{I_F})$. Since the latter is a subgroup of $\widehat{k^*_F}$, it is finite and abelian.   \\

Let $g\in C_{PGL_N(\mathbb{C})}(\rho|_{I_F},\rho(Fr)C^0_{PGL_N(\mathbb{C})}(\rho|_{I_F}),E)$, and let $\tilde{g}$ be a representative of $g$ in  $GL_N(\mathbb{C})$. We set
\begin{equation}\label{chig}
    \chi_g:=\tilde{g}(\tilde{\rho}|_{I_F})\tilde{g}^{-1}(\tilde{\rho}^{-1}|_{I_F}): I_F\rightarrow GL_N(\mathbb{C})
\end{equation}
We observe that $\chi_g$ does not depend on the choice of $\tilde{g}$.

Since  in particular $g\in C_{PGL_N(\mathbb{C})}(\rho|_{I_F})$, there holds
\[\eta(\chi_g)(\omega)=\eta(\tilde{g}\tilde{\rho}(\omega)\tilde{g}^{-1}\tilde{\rho}(\omega)^{-1})=1\]
for any $\omega\in I_F$, so $\chi_g(\omega)$ is a scalar for any $\omega\in I_F$. In this way we have a character
\[\chi_g:I_F\rightarrow \mathbb{C}^*.\]
We show that it is Frobenius stable.
Since $g\in C_{PGL_N(\mathbb{C})}(\rho(Fr)C^0_{PGL_N(\mathbb{C})}(\rho|_{I_F}))$, we have 
 \[c:=\rho(Fr)^{-1}g\rho(Fr)g^{-1}\in C^0_{PGL_N(\mathbb{C})}(\rho|_{I_F}).\]
 Lifting this identity to $GL_N(\mathbb{C})$ yields
 \begin{equation}\label{commid}
 \tilde{c}:=\tilde{\rho}(Fr)^{-1}\tilde{g}\tilde{\rho}(Fr)\tilde{g}^{-1}\in C_{GL_n(\mathbb{C})}(\tilde{\rho}|_{I_F})
 \end{equation}
 Therefore, for any $\omega\in I_F$ 
 \begin{align*}
  &\chi_g(Fr \omega Fr^{-1})=& \\
    &\tilde{g}(\tilde{\rho}(Fr \omega Fr^{-1}))\tilde{g}^{-1}(\tilde{\rho}(Fr \omega Fr^{-1})^{-1})=&\\
    &\tilde{g}\tilde{\rho}(Fr)(\tilde{\rho}( \omega ))\tilde{\rho}(Fr^{-1})\tilde{g}^{-1}\tilde{\rho}^{-1}(Fr^{-1})(\tilde{\rho}( \omega )^{-1})\tilde{\rho}^{-1}(Fr)=\\
    &\tilde{g}\tilde{\rho}(Fr)
    (\tilde{\rho}( \omega ))
    \tilde{\rho}(Fr)^{-1}\tilde{g}^{-1}
    \tilde{\rho}(Fr)(\tilde{\rho}( \omega )^{-1})\tilde{\rho}(Fr)^{-1}= & \tilde{g}\tilde{\rho}(Fr)=\tilde{\rho}(Fr)\tilde{c}\tilde{g} \text{ \  by \eqref{commid}}\\
    &\tilde{\rho}(Fr)\tilde{c}\tilde{g}
    (\tilde{\rho}( \omega ))
    \tilde{g}^{-1}\tilde{c}^{-1}\tilde{\rho}(Fr)^{-1}
    \tilde{\rho}(Fr)(\tilde{\rho}( \omega )^{-1})\tilde{\rho}(Fr)^{-1}=\\
     &\tilde{\rho}(Fr)\tilde{c}\tilde{g}
    (\tilde{\rho}( \omega ))
    \tilde{g}^{-1}\tilde{c}^{-1}
    (\tilde{\rho}( \omega )^{-1})\tilde{\rho}(Fr)^{-1}= &\tilde{c}\in C_{GL_N(\mathbb{C})}(\tilde{\rho})\text{ \  by \eqref{commid}}\\
     &\tilde{\rho}(Fr)\tilde{c}\tilde{g}
    (\tilde{\rho}( \omega ))
    \tilde{g}^{-1}(\tilde{\rho}( \omega )^{-1})\tilde{c}^{-1}\tilde{\rho}(Fr)^{-1}=&\\
    &\tilde{\rho}(Fr)\tilde{c}\chi_g( \omega )
    \tilde{c}^{-1}\tilde{\rho}(Fr)^{-1}=\chi_g(\omega).   &\chi_g(\omega)\in \mathbb{C}^* 
 \end{align*}

Hence $\chi_g$ is a Frobenius stable character of $I_F$. Moreover by definition $\chi_g|_{P_F}=1$, hence we can regard $\chi_g$ as a character of $k_F^*$. By construction of $\chi_g$ it holds 
\[\chi_g\tilde{\rho}|_{I_F}=\tilde{g}(\tilde{\rho}|_{I_F})\tilde{g}^{-1}.\]
In addition $g\in C_{PGL_N(\mathbb{C})}(E)$ and hence $\tilde{g}\in C_{GL_N(\mathbb{C})}(E)$, so $(\tilde{\rho},E)\sim_{I_F}(\tilde{\rho}\otimes\chi_g,E)$ in $(\Phi_N)_0$. 

Therefore the assignment $g\mapsto\chi_g$ gives a map
\begin{align}\label{chi}
    \widetilde{\Xi}: C_{PGL_N(\mathbb{C})}(\rho|_{I_F},\rho(Fr)C^0_{PGL_N(\mathbb{C})}(\rho|_{I_F}),E)&\rightarrow Stab_{\widehat{k^*_F}}((\tilde{\rho},E)_{I_F}).
\end{align}
The map $\widetilde{\Xi}$ is a group morphism: for any $g_1,g_2\in C_{PGL_N(\mathbb{C})}(\rho|_{I_F},\rho(Fr)C^0_{PGL_N(\mathbb{C})}(\rho|_{I_F}),E)$ 
\begin{align*}
    \chi_{g_1g_2}\tilde{\rho}|_{I_F}=\tilde{g_1}\tilde{g_2}(\tilde{\rho}|_{I_F})\tilde{g_2}^{-1}\tilde{g_1}^{-1}=\tilde{g_1}(\chi_{g_2} \tilde{\rho}|_{I_F})\tilde{g_1}^{-1}=\chi_{g_2}\tilde{g_1}(\tilde{\rho}|_{I_F})\tilde{g_1}^{-1}=\chi_{g_1}\chi_{g_2}\tilde{\rho}|_{I_F}
\end{align*}
hence $\chi_{g_1g_2}=\chi_{g_1}\chi_{g_2}$.

It is also surjective. Indeed let $\chi\in \widehat{k^*_F}$  such that $(\chi\otimes\tilde{\rho},E)_{I_F}=(\tilde{\rho},E)_{I_F}$. There exists an element $\tilde{A}\in GL_N(\mathbb{C})$ satisfying conditions \eqref{aGln}, \eqref{b} and \eqref{c} of Definition \ref{IFequivsln}.
This implies that $A=\eta(\tilde{A})\in C_{PGL_N(\mathbb{C})}(\rho|_{I_F}, \rho(Fr)C^0_{PGL_N(\mathbb{C})}(\rho|_{I_F}),E)$. From $\tilde{A}(\tilde{\rho}|_{I_F})\tilde{A}^{-1}=\chi\otimes(\tilde{\rho}|_{I_F})$ we obtain 
\[\chi=\tilde{A}(\tilde{\rho}|_{I_F})\tilde{A}^{-1} (\tilde{\rho}|_{I_F})^{-1}=\chi_{A}=\widetilde{\Xi}(A).\]

Since the image of $\widetilde{\Xi}$ is a finite group, the morphism \eqref{chi} factors through the component group
\begin{equation}\label{chibar}
\begin{tikzcd}
{C_{PGL_N(\mathbb{C})}(\rho|_{I_F},\rho(Fr)C^0_{PGL_N(\mathbb{C})}(\rho|_{I_F}),E)} \arrow[d] \arrow[rd, "\widetilde{\Xi}"]                                                                               &                                                                                  \\
{A_{PGL_N(\mathbb{C})}(\rho|_{I_F},\rho(Fr)C^0_{PGL_N(\mathbb{C})}(\rho|_{I_F}),E)} \arrow[r, "\Xi"] & Stab_{\widehat{k^*_F}}((\tilde{\rho},E)_{I_F})
\end{tikzcd}
\end{equation}
\nomenclature{$\Xi$}{The isomorphism in \eqref{chibar}}
So the map $\Xi$ in diagram \eqref{chibar} is a surjective group morphism. We prove that it is also injective by showing that the $ker(\tilde{\Xi})=C^0_{PGL_N(\mathbb{C})}(\rho|_{I_F},\rho(Fr)C^0_{PGL_N(\mathbb{C})}(\rho|_{I_F}),E)$. 

Let $g\in ker(\tilde{\Xi})$ and let $\tilde{g}$ be a representative of $g$ in $GL_N(\mathbb{C})$. Then 
\[\tilde{g}(\tilde{\rho}|_{I_F})\tilde{g}^{-1}(\tilde{\rho}^{-1}|_{I_F})=1,\] that is $\tilde{g}\in C_{GL_N({\mathbb{C})}}(\tilde{\rho}|_{I_F})$. Since $g\in C_{PGL_N(\mathbb{C})}(E)$, it holds $\tilde{g}\in C_{GL_N(\mathbb{C})}(\tilde{\rho}|_{I_F},E)$. The latter is connected because it is the centralizer of an element in $GL_N(\mathbb{C})$.

\noindent By Remark \ref{rmkIFGln},  for any $\tilde{g}\in C_{GL_N(\mathbb{C})}(\tilde{\rho}|_{I_F},E)$ it holds $\tilde{g}\tilde{\rho}(Fr) C_{GL_N(\mathbb{C})}(\tilde{\rho}|_{I_F})\tilde{g}^{-1}=\tilde{\rho}(Fr) C_{GL_N(\mathbb{C})}(\tilde{\rho}|_{I_F})$, so 
\[\tilde{g}\in C_{GL_N(\mathbb{C})}(\rho|_{I_F},E)= C_{GL_N(\mathbb{C})}(\tilde{\rho}|_{I_F},\tilde{\rho}(Fr)C_{GL_N(\mathbb{C})}(\tilde{\rho}|_{I_F}),E).\]
Then 
\[g=\eta(\tilde{g})\in \eta(C_{GL_N(\mathbb{C})}(\tilde{\rho}|_{I_F},\tilde{\rho}(Fr)C_{GL_N(\mathbb{C})}(\tilde{\rho}|_{I_F}),E))\leq C_{PGL_N(\mathbb{C})}(\rho|_{I_F},\rho(Fr)C^0_{PGL_N(\mathbb{C})}(\rho|_{I_F}),E).\]
\sloppy Since $C_{GL_N(\mathbb{C})}(\tilde{\rho}|_{I_F},\tilde{\rho}(Fr)C_{GL_N(\mathbb{C})}(\tilde{\rho}|_{I_F}),E)$ is the same as $ C_{GL_N(\mathbb{C})}(\tilde{\rho}|_{I_F},E)$, it is connected. It follows that 
\[\eta(C_{GL_N(\mathbb{C})}(\tilde{\rho}|_{I_F},\tilde{\rho}(Fr)C_{GL_N(\mathbb{C})}(\tilde{\rho}|_{I_F}),E))\leq C^0_{PGL_N(\mathbb{C})}(\rho|_{I_F},\rho(Fr)C^0_{PGL_N(\mathbb{C})}(\rho|_{I_F}),E)\]
that is,
\[ g\in C^0_{PGL_N(\mathbb{C})}(\rho|_{I_F},\rho(Fr)C^0_{PGL_N(\mathbb{C})}(\rho|_{I_F}),E),\]
giving injectivity of $\Xi$.

\end{proof}
\end{theorem}
\subsection{Comparison of the fibers of $\mathcal{M}'_N$ and $\mathcal{L}'_N$}
The aim of this section is to show that the fibers of the Langlands correspondence and the fibers of the Macdonald-Vogan correspondence may not have the same cardinality, and the natural map \eqref{iotahat} between them is neither an injection nor a surjection in general.

Let $(\rho,E)\in (\Phi'_N)_0$. By \cite[Theorem 4.3]{Gel} the component group $A_{PGL_N(\mathbb{C})}(\rho, E)$ is finite abelian and there is a canonical simply transitive action of its character group on ${\mathcal{L}'_N}^{-1}(\rho,E)$. Therefore there is a bijection between ${\mathcal{L}'_N}^{-1}(\rho,E)$ and the group $A_{PGL_N(\mathbb{C})}(\rho, E)^{\wedge}$.

\sloppy Similarly, by Proposition \ref{fiber2} there is a bijection between ${\mathcal{M}'_N}^{-1}((\rho,E)_{I_F})$ and the group
$A_{PGL_N(\mathbb{C})}(\rho|_{I_F},\rho(Fr)C^0_{PGL_N(\mathbb{C})}(\rho|_{I_F}),E)^{\wedge}$.

The inclusion $C_{PGL_N(\mathbb{C})}(\rho, E)\hookrightarrow C_{PGL_N(\mathbb{C})}(\rho|_{I_F},\rho(Fr)C^0_{PGL_N(\mathbb{C})}(\rho|_{I_F}),E)$ induces a  map between the component groups
\begin{equation}\label{inclusioncomp}
\iota:A_{PGL_N(\mathbb{C})}(\rho, E)\rightarrow A_{PGL_N(\mathbb{C})}(\rho|_{I_F},\rho(Fr)C^0_{PGL_N(\mathbb{C})}(\rho|_{I_F}),E)
\end{equation}
\nomenclature{$\iota$}{The map in \eqref{inclusioncomp}}
that gives by duality a map
\begin{equation}\label{iotahat}
    \hat{\iota}: A_{PGL_N(\mathbb{C})}(\rho|_{I_F},\rho(Fr)C^0_{PGL_N(\mathbb{C})}(\rho|_{I_F}),E)^{\wedge}\rightarrow A_{PGL_N(\mathbb{C})}(\rho, E)^{\wedge}.
\end{equation}
This map is neither injective nor surjective in general. We give two examples to illustrate this.

We will denote by $T$ the diagonal torus in $PGL_N(\mathbb{C})$. 

\begin{enumerate}
\item \textbf{Failure of the surjectivity}: 
Let $(\rho,E)\in(\Phi'_N)_0$ be such that
\begin{itemize}
\item $\rho|_{I_F}=1$
\item $\rho(Fr)=\begin{pmatrix}1 &&&\\ &\zeta&&\\&&\ddots&\\&&&\zeta^{N-1}\end{pmatrix}$ where $\zeta\in\mathbb{C}$ a primitive $N$-th root of $1$
\item $E=0$
\end{itemize}
In this case
\[A_{PGL_N(\mathbb{C})}(\rho, E)=\faktor{C_{PGL_N(\mathbb{C})}(\rho(Fr))}{C^0_{PGL_N(\mathbb{C})}(\rho(Fr))}\]
and it holds
\[C_{PGL_N(\mathbb{C})}(\rho(Fr))=T\rtimes\left\langle\begin{psmallmatrix}0&1&&&\\&0&1&&\\&&\ddots&\ddots&\\0&\dots &\dots  &0&1\\
1&0&\dots&\dots &0\end{psmallmatrix}\right\rangle.\]
Therefore $A_{PGL_N(\mathbb{C})}(\rho, E)$ is the cyclic group of order $N$.
So  $|{\mathcal{L}'_N}^{-1}(\rho,E)|=N$.

On the other hand, in this case $C^0_{PGL_N(\mathbb{C})}(\rho|_{I_F})=C^0_{PGL_N(\mathbb{C})}(1)=PGL_N(\mathbb{C})$, hence 
\begin{align*}
    C_{PGL_N(\mathbb{C})}(\rho|_{I_F},\rho(Fr)C^0_{PGL_N(\mathbb{C})}(\rho|_{I_F}),E)&=C_{PGL_N(\mathbb{C})}(1,PGL_N(\mathbb{C}),0)\\
    &=PGL_N(\mathbb{C})=C^0_{PGL_N(\mathbb{C})}(\rho|_{I_F},\rho(Fr)C^0_{PGL_N(\mathbb{C})}(\rho|_{I_F}),E)
\end{align*}
 
Hence $|{\mathcal{M}'_N}^{-1}((\rho,E)_{I_F})|=|A_{PGL_N(\mathbb{C})}(\rho|_{I_F},\rho(Fr)C^0_{PGL_N(\mathbb{C})}(\rho|_{I_F}),E)^{\wedge}|=1$.

In this case, the map
\[\hat{\iota}: 1\cong  A_{PGL_n(\mathbb{C})}(\rho|_{I_F},\rho(Fr)C^0_{PGL_n(\mathbb{C})}(\rho|_{I_F}),E)^{\wedge}\rightarrow  A_{PGL_N(\mathbb{C})}(\rho, E)^{\wedge}\]
is  injective but not surjective.

For the sake of completeness, we describe explicitly the representations in ${\mathcal{M}'_N}^{-1}((\rho,E)_{I_F})$ and ${\mathcal{L}'_N}^{-1}(\rho,E)$.
 Let $\chi_{\zeta}\in \widehat{F^*}$ be the unramified character of $F^*$ that takes value $\zeta$ on any uniformizer of $\widehat{F^*}$. This corresponds by the bijection of local class field theory to the unramified character of $W_F$ that takes value $\zeta$ at $Fr$. Let $\tilde{T}$ be the diagonal torus in $G_N$. Then $\tilde{\pi}:=\pind_{\tilde{T}}^{G_N}(\bigboxtimes_{i=0}^{n-1} \chi_{\zeta}^{i})$ is an irreducible representation of $G_N$. The representations in ${\mathcal{L}'_N}^{-1}(\rho,E)$ are the $N$  irreducible constituents of $Res^{G_N}_{G'_N}\tilde{\pi}$.

 On the other hand, in this case ${\mathcal{M}'_N}^{-1}((\rho,E)_{I_F})$ contains only the trivial representation.

\item\textbf{Failure of injectivity:}
Assume that $(q-1,N)\neq 1$. Let $e>1$ be a common divisor of $q-1$ and $N$ and let $\zeta\in\mathbb{C}$ be a primitive $e$-th root of $1$.


We consider $(\rho, E)\in (\Phi'_N)_0$ such that
\begin{itemize}
\item  $\rho(I_F)=\langle M\rangle$ with 
\[M:=\begin{pmatrix}\mathbf{1}_{\frac{N}{e}}&&&\\&\zeta \mathbf{1}_{\frac{N}{e}}&&\\&&\ddots&\\&&&\zeta^{e-1}\mathbf{1}_{\frac{N}{e}}\end{pmatrix}\]
where $\mathbf{1}_{\frac{N}{e}}$ denotes the identity matrix of dimension $\frac{N}{e}$;
\item $\rho(Fr)$ is a regular element of $T$
\item $E=0$.
\end{itemize}
By construction, $M^q=M$, so a Weil-Deligne representation just needs  to map $Fr$ to an element commuting with $M$, and $\rho(Fr)$ satisfies this condition.
We have 
\[C_{PGL_N(\mathbb{C})}(M)\cong \faktor{(GL_{\frac{N}{e}}(\mathbb{C})^e\rtimes \faktor{\mathbb{Z}}{e\mathbb{Z}})}{\mathbb{C}^*}.\] Moreover $C_{PGL_N(\mathbb{C})}(\rho(Fr))=T$.
 So we have 
\[C_{PGL_N(\mathbb{C})}(\rho, E)= C_{PGL_N(\mathbb{C})}(M, \rho(Fr))= T=C^0_{PGL_N(\mathbb{C})}(\rho, E).\]
Hence 
\[A_{PGL_N(\mathbb{C})}(\rho, E)=1.\]
On the other hand,
\[C^0_{PGL_N(\mathbb{C})}(\rho|_{I_F})=C^0_{PGL_N(\mathbb{C})}(M)=\faktor{GL_{\frac{N}{e}}(\mathbb{C})^e}{\mathbb{C}^*}.\] 
This group contains $T$, and hence $\rho(Fr)$, so we have
\[C_{PGL_N(\mathbb{C})}(\rho|_{I_F},\rho(Fr)C^0_{PGL_N(\mathbb{C})}(\rho|_{I_F}),E)=C_{PGL_N(\mathbb{C})}(M,C^0_{PGL_N(\mathbb{C})}(M))=C_{PGL_N(\mathbb{C})}(M).\]
Hence,
\[A_{PGL_N(\mathbb{C})}(\rho|_{I_F},\rho(Fr)C^0_{PGL_N(\mathbb{C})}(\rho|_{I_F}),E)=\faktor{C_{PGL_N(\mathbb{C})}(M)}{C^0_{PGL_N(\mathbb{C})}(M)}\cong \faktor{\mathbb{Z}}{e\mathbb{Z}}\]

    Since we assumed $e>1$, it follows that in this case, the map
\[\hat{\iota}:A_{PGL_n(\mathbb{C})}(\rho|_{I_F},\rho(Fr)C^0_{PGL_n(\mathbb{C})}(\rho|_{I_F}),E)^{\wedge} \rightarrow  A_{PGL_N(\mathbb{C})}(\rho, E)^{\wedge}\cong 1\]
is surjective but not injective.

As in the previous example, we describe explicitly the representations in ${\mathcal{M}'_N}^{-1}((\rho,E)_{I_F})$ and ${\mathcal{L}'_N}^{-1}(\rho,E)$. For the sake of simplicity of the description, assume $\rho(Fr)$ to be an element of the compact torus, i.e. assume the entries of $\rho(Fr)$ to have all the same absolute value in $\mathbb{C}$.

Let $\tilde{\rho}$ be a lift of $\rho$ to $GL_N(\mathbb{C})$ such that $ \tilde{\rho}(Fr)=\begin{psmallmatrix} {t_1}&&&\\ &{t_2}&&\\&&\ddots&\\&&&{t_N}\end{psmallmatrix}$ with $|t_j|=1$ for any $j=1\dots N$. Then $\tilde{\rho}=\bigoplus_{j=1}^N \chi_j$ with $\chi_j\in \widehat{W_F}$. By the bijection of local class field theory, any $\chi_j$ can be regarded as a tame character of $F^*$. More explicitly, $\chi_{j}$ takes value $t_j$ on any uniformizer of $\widehat{F^*}$, and $\chi_j|_{\of{F}^*}=Infl_{\faktor{\of{F}^*}{1+\pf{F}}}^{\of{F}^*} \widehat{\zeta}^{\lfloor\frac{e(j-1)}{N} \rfloor}$, where $\widehat{\zeta}$ denotes a character of $k_F^*\cong \faktor{\of{F}^*}{1+\pf{F}}$ of order $e$.
Then $\tilde{\pi}:=\pind_{G_1^N}^{G_N}(\bigboxtimes_{i=0}^{n-1} \chi_{j})$ is a tempered irreducible representation of $G_N$. This representation remains irreducible upon restriction to $G'_N$, and   $Res^{G_N}_{G'_N}\tilde{\pi}$ is the unique element in ${\mathcal{L}'_N}^{-1}(\rho,E)$. \\

 On the other hand, the $I_F$-equivalence class of $(\tilde{\rho},0)$ corresponds, according to Macdonald's construction for $G_N$, to the partition valued function on $\overline{\mathcal{C}}$ taking value $(1^{\frac{N}{e}})$ on $\widehat{\zeta}^i$ for any $0\leq i\leq e-1$, and $0$ everywhere else. Hence ${\mathcal{M}_N}^{-1}((\tilde{\rho},0)_{I_F})=\pind_{G_{\frac{N}{e}}^e}^{G_N}(\bigboxtimes_{i=0}^{e-1} (\widehat{\zeta}^i\circ det))$. This is an irreducible representation of $G_N$, but when restricted to ${G'_N}$ it splits in $e$ inequivalent irreducible subrepresentations. These are exactly the representations into ${\mathcal{M}'_N}^{-1}((\rho,E)_{I_F})$.
\end{enumerate}

\section{Compatibility between Langlands and Macdonald-Vogan parameterizations}
\subsection{Compatibility between $\mathcal{L}'_N$ and $\mathcal{M}'_N$}
The main result of this section is Theorem \ref{compatibilitysln}, which yields compatibility of Langlands and Macdonald correspondences, considering all the representations in the same fiber together.
\subsubsection{The head of Parahoric Restriction}

 We recall some notation from Section \ref{Gln}. Let $\tilde{\pi}\in (\Omega_N)_0$ be a representation of $G_N$. The bijection \eqref{deflambdapi} associates to $\tilde{\pi}$ the partition valued function $\Lambda_{\tilde{\pi}}:\mathcal{C}_0\rightarrow Par$ of degree $N$, and the map \eqref{deflambdabar} associates to $\Lambda_{\tilde{\pi}}$ a partition valued function $\overline{\Lambda_{\tilde{\pi}}}:\overline{\mathcal{C}}\rightarrow Par$ of degree $N$. The bijection \eqref{sigmabar} associates to $\overline{\Lambda_{\tilde{\pi}}}$  the irreducible representation  $\overline{\sigma}_{\overline{\Lambda_{\tilde{\pi}}}}$ of $\overline{G}_N$. 
By Theorem \ref{prephead} the irreducible constituents of $\mathcal{P}^{G_N}_{\overline{G}_N}(\tilde{\pi})$ are of the form
$\overline{\sigma}_{\overline{\Lambda}}$ for $\Lambda\leq\Lambda_{\tilde{\pi}}$, and $\mathcal{HP}^{G_N}_{\overline{G}_N}(\tilde{\pi})=\overline{\sigma}_{\overline{\Lambda_{\tilde{\pi}}}}$ has multiplicity $1$. 
 
 Let $\chi\in\widehat{F^*}$. Since $\chi$ is tame, the restriction $Res^{F^*}_{\of{F}^*}\chi$  gives a character $\overline{\chi}\in \widehat{\faktor{\of{F}^*}{1+\pf{F}}}\cong \widehat{k_F^*}$. The character group $\widehat{k_F^*}$ acts on the set of partition valued functions on $\overline{\mathcal{C}}$  by 
 \begin{equation}\label{actiononlambda}
 (\overline{\chi}\overline{\Lambda})(\overline{\tau}):=\overline{\Lambda}( \overline{\tau}\otimes(\overline{\chi}^{-1}\circ det)).
 \end{equation}
 
\begin{lemma}\label{twistedparres}
  With the above notation, the irreducible constituents of  $\mathcal{P}^{G_N}_{\overline{G}_N}(\tilde{\pi}\otimes (\chi\circ\det))$ are of the form  $\overline{\sigma}_{\overline{\chi}\overline{\Lambda}}$ with $\Lambda\leq\Lambda_{\tilde{\pi}}$.  
 \begin{proof}
 As a first step, we show that $\mathcal{P}^{G_N}_{\overline{G}_N}(\tilde{\pi}\otimes (\chi\circ\det))=\mathcal{P}^{G_N}_{\overline{G}_N}(\tilde{\pi})\otimes (\overline{\chi}\circ\det)$. Indeed since $\chi|_{1+\pf{F}}=1$, there holds $\chi\circ det|_{K_N^+}=1$, therefore $(\tilde{\pi}\otimes (\chi\circ\det)) ^{K_N^+}=\tilde{\pi}^{K_N^+}\otimes (\chi\circ\det)$. Moreover $det(K_N)\subseteq \of{F}^*$, so
 \begin{align*}
     \mathcal{P}^{G_N}_{\overline{G}_N}(\tilde{\pi}\otimes (\chi\circ\det)) =&Res^{G_N}_{K_N}(\tilde{\pi}\otimes (\chi\circ\det)) ^{K_N^+} =\\
     &Res^{G_N}_{K_N}(\tilde{\pi})^{K_N^+}\otimes  Res^{G_N}_{K_N}(\chi\circ\det)^{K_N^+} = \\ &Res^{G_N}_{K_N}(\tilde{\pi})^{K_N^+}\otimes (\overline{\chi}\circ\det) = .  \mathcal{P}^{G_N}_{\overline{G}_N}(\tilde{\pi})\otimes (\chi\circ\det)
     \end{align*}
  
 Therefore the irreducible constituents of  
 $\mathcal{P}^{G_N}_{\overline{G}_N}(\tilde{\pi}\otimes (\chi\circ\det))$ are of the form $\overline{\sigma}_{\overline{\Lambda}}\otimes (\overline{\chi}\circ det)$ for $\Lambda\leq\Lambda_{\tilde{\pi}}$.\\
 
 Next we show that for any $\Lambda\leq\Lambda_{\overline{\sigma}}$ and $\overline{\chi}\in\widehat{k_F^*}$
 \begin{equation}\label{Hpdet}
 \overline{\sigma}_{\overline{\Lambda}}\otimes (\overline{\chi}\circ det)= \overline{\sigma}_{\overline{\chi}\overline{\Lambda}}.
 \end{equation}
 The representations $\overline{\sigma}_{\overline{\Lambda}}$ are defined in \eqref{sigmabar} as $\overline{\sigma}_{\overline{\Lambda}}=\pind_{\prod_{\overline{\tau}\in\overline{\mathcal{C}}}\overline{G}_{d(\overline{\tau})|\overline{\Lambda}(\overline{\tau})|}}^{\overline{G}_N}(\bigboxtimes_{\overline{\tau}\in\overline{\mathcal{C}}}\overline{\tau}^{\overline{\Lambda}(\overline{\tau})})$, where $\overline{\tau}^{\overline{\Lambda}(\overline{\tau})}$ is the irreducible constituent of $\pind_{\overline{G}_{d(\overline{\tau})}^{|\overline{\Lambda}(\overline{\tau})|}}^{\overline{G}_{d(\overline{\tau})|\overline{\Lambda}(\overline{\tau})|}}(\bigboxtimes_{i=1}^{|\overline{\Lambda}(\overline{\tau})|}\overline{\tau})$ corresponding to the partition $\overline{\Lambda}(\overline{\tau})$.
It follows that
 \begin{equation}\label{const1}
     \overline{\tau}_{\overline{\Lambda}}\otimes (\overline{\chi}\circ det)=\pind_{\prod_{\overline{\tau}\in\overline{\mathcal{C}}}\overline{G}_{d(\overline{\tau})|\overline{\Lambda}(\overline{\tau})|}}^{\overline{G}_N}(\bigboxtimes_{\overline{\tau}\in\overline{\mathcal{C}}}\overline{\tau}^{\overline{\Lambda}(\overline{\tau})})\otimes(\overline{\chi}\circ\det)=
\pind_{\prod_{\overline{\tau}\in\overline{\mathcal{C}}}\overline{G}_{d(\overline{\tau})|\overline{\Lambda}(\overline{\tau})|}}^{\overline{G}_N}(\bigboxtimes_{\overline{\tau}\in\overline{\mathcal{C}}}(\overline{\tau}^{\overline{\Lambda}(\overline{\tau})}\otimes(\overline{\chi}\circ\det))).
 \end{equation}
On the other hand,
\begin{equation}\label{const2}
\overline{\sigma}_{\overline{\chi}\overline{\Lambda}}=
\pind_{\prod_{\overline{\tau}\in\overline{\mathcal{C}}}\overline{G}_{d(\overline{\tau})|\overline{\chi}\overline{\Lambda}(\overline{\tau})|}}^{\overline{G}_N}(\bigboxtimes_{\overline{\tau}\in\overline{\mathcal{C}}}\overline{\tau}^{\overline{\chi}\overline{\Lambda}(\overline{\tau})})=
\pind_{\prod_{\overline{\tau}\in\overline{\mathcal{C}}}\overline{G}_{d(\overline{\tau})|\overline{\Lambda}(\overline{\tau})|}}^{\overline{G}_N}(\bigboxtimes_{\overline{\tau}\in\overline{\mathcal{C}}}(\overline{\tau}\otimes(\overline{\chi}\circ\det))^{\overline{\Lambda}(\overline{\tau})}).
\end{equation}
 We show that for any $\overline{\tau}\in\overline{\mathcal{C}}$ and any $\Lambda\leq\Lambda_{\tilde{\pi}}$ it holds
\begin{equation}\label{twist}
\overline{\tau}^{\overline{\Lambda}(\overline{\tau})}\otimes(\overline{\chi}\circ\det)=(\overline{\tau}\otimes(\overline{\chi}\circ\det))^{\overline{\Lambda}(\overline{\tau})}.
\end{equation}

Let $d:=d(\overline{\tau})$ and $n:=|\overline{\Lambda}(\overline{\tau})|$. We have 
\[\pind_{\overline{G}_{d}^{n}}^{\overline{G}_{dn}}(\bigboxtimes_{i=1}^{n}\overline{\tau})\otimes(\overline{\chi}\circ\det)=\pind_{\overline{G}_{d}^{n}}^{\overline{G}_{dn}}(\bigboxtimes_{i=1}^{n}(\overline{\tau}\otimes(\overline{\chi}\circ\det))).\]
It implies that the 
$End_{\overline{G}_{dn}}\bigg(\pind_{\overline{G}_{d}^{n}}^{\overline{G}_{dn}}\big(\bigboxtimes_{i=1}^{n}(\overline{\tau}\otimes(\overline{\chi}\circ\det))\big)\bigg)$
-module structure over \\
$Hom_{\overline{G}_{dn}}\bigg(\pind_{\overline{G}_{d}^{n}}^{\overline{G}_{dn}}\big(\bigboxtimes_{i=1}^{n}(\overline{\tau}\otimes(\overline{\chi}\circ\det))\big),
\overline{\tau}^{\overline{\Lambda}(\tau)}\otimes(\overline{\chi}\circ\det)\bigg)$
is the same as the 
$End_{\overline{G}_{dn}}\big(\pind_{\overline{G}_{d}^{|\overline{\Lambda(\overline{\tau})|}}}^{\overline{G}_{dn}}(\bigboxtimes_{i=1}^{n}\overline{\tau})\bigg)$
-module structure over 
$Hom_{\overline{G}_{dn}}\bigg(\pind_{\overline{G}_{d}^{|\overline{\Lambda}(\overline{\tau})|}}^{\overline{G}_{dn}}\big(\bigboxtimes_{i=1}^{n}\overline{\tau}\big),\overline{\tau}^{\Lambda(\overline{\tau})}\bigg)$.
This proves \eqref{twist}.  Comparing \eqref{const1} and \eqref{const2}, it follows that $ \overline{\sigma}_{\overline{\Lambda}}\otimes (\overline{\chi}\circ det)= \overline{\sigma}_{\overline{\chi}\overline{\Lambda}}$, hence the statement.
  \end{proof}
\end{lemma}

\begin{cor}\label{twistedparrescor}
For any  representation $\tilde{\pi}\in(\Omega_N)_0$ and for any  tame  character $\chi\in\widehat{F^*}$, it holds
\[\mathcal{HP}^{G_N}_{\overline{G}_N}(\tilde{\pi}\otimes(\chi\circ \det)) = \mathcal{HP}^{G_N}_{\overline{G}_N}(\tilde{\pi})\otimes(\overline{\chi}\circ \det)\]
where $\overline{\chi}$ is the character of $k_F^*$ obtained by restriction of $\chi$ to $\of{F}^*$.
\begin{proof}
We retain the notation from Lemma \ref{twistedparres}. The description of the irreducible constituents of  $\mathcal{P}^{G_N}_{\overline{G}_N}(\tilde{\pi}\otimes\chi\circ \det) $ gives 
    \[\mathcal{HP}^{G_N}_{\overline{G}_N}(\tilde{\pi}\otimes(\chi\circ \det))=\overline{\sigma}_{\overline{\chi}\overline{\Lambda_{\tilde{\pi}}}}.\]
    and by \eqref{Hpdet}
    \[\overline{\sigma}_{\overline{\chi}\overline{\Lambda_{\tilde{\pi}}}}=
     \overline{\sigma}_{\overline{\Lambda}_{\tilde{\pi}}}\otimes (\overline{\chi}\circ det)=\mathcal{HP}^{G_N}_{\overline{G}_N}(\tilde{\pi})\otimes(\chi\circ \det).\]
\end{proof}
\end{cor}
For any $\pi\in (\Omega'_N)_0$, we denote by $F^*\pi$ the orbit of $\pi$ under the conjugation action of the group $\faktor{G_N}{G'_N}\cong F^*$. This orbit is a finite set \cite[Theorem 4.1]{Gel}. We set
\begin{equation}\label{oplus}
    \pi^\oplus:=\bigoplus_{\pi'\in F^*\pi}\pi'.
\end{equation}
\nomenclature{$\pi^{\oplus}$}{Sum of all the representation in the same $L$-packet as $\pi$ \eqref{oplus}}

Similarly, for any $\overline\pi\in \overline{\Omega_N'}$, we denote by $k_F^*\overline\pi$ the orbit of $\overline\pi$ under the conjugation action of the group $\faktor{\overline{G}_N}{\overline{G'}_N}\cong k_F^*$, and we set 
\begin{equation}\label{oplusbar}
    \overline\pi^\oplus:=\bigoplus_{\overline\pi'\in k_F^*\overline\pi}\overline\pi'.
\end{equation}
\nomenclature{$\overline{\pi}^{\oplus}$}{Sum of all the representation in the same Macdonal-Vogan fiber as $\overline{\pi}$ \eqref{oplusbar}}
Since the restriction of irreducible representations from $G_N$ (respectively $\overline{G}_N$) to $G'_N$ (respectively $\overline{G'}_N$ ) is multiplicity free, it holds
\begin{align}\label{plusrestrictiongln}
    &\pi^\oplus:=Res^{G_N}_{G'_N} R(\pi) &\overline\pi^\oplus:=Res^{\overline{G}_N}_{\overline{G'}_N}\overline{R}(\overline\pi)
\end{align}
where $R(\pi)$ (respectively $\overline{R}(\overline{\pi})$) is an irreducible representation of $G_N$ (respectively of $\overline{G}_N$) containing  $\pi$ (respectively $\overline{\pi})$ as $G'_N$-subrepresentation (respectively $\overline{G'_N}$-subrepresentation).

Moreover by definition of the maps $\mathcal{L}'_N$ and $\mathcal{M}'_N$, for any $\pi\in(\Omega'_N)_0$ and for any $\overline{\pi}\in \overline{\Omega'_N}$ it holds
\begin{align}\label{plusmaclang}
&\pi^\oplus =\bigoplus_{\pi'\in {\mathcal{L}'_N}^{-1}(\mathcal{L}'_N(\pi))}\pi' &\text{and} \hspace{1cm} &\overline\pi^\oplus =\bigoplus_{\overline\pi'\in {\mathcal{M}'_N}^{-1}(\mathcal{M}'_N(\overline\pi))}\overline\pi'
\end{align}

 We extend the notion of head of parahoric restriction from $G_N$ to representations of $G'_N$ of the form $\pi^{\oplus}$ for some $\pi\in(\Omega'_N)_0$ as follows.

 \begin{defi}\label{HPSLn}
     Let $\pi\in(\Omega'_N)_0$ , let $R(\pi)\in(\Omega_N)_0$ be a representation of $G_N$ containing $\pi$ as $G'_N$-subrepresentation. The head of the parahoric restriction of $\pi^{\oplus}$ is the  $\overline{G'}_N$-representation  
\[\mathcal{HP}^{G'_N}_{\overline{G'}_N}\pi^{\oplus}:=Res^{\overline{G}_N}_{\overline{G'}_N}\mathcal{HP}^{G_N}_{\overline{G}_N}R(\pi)\]
 \end{defi}

The definition above is well posed: by Corollary \ref{twistedparrescor}, the representation $\mathcal{HP}^{G'_N}_{\overline{G'}_N}\pi^{\oplus}$ does not depend on the choice of $R(\pi)$ .

\begin{lemma}\label{HPsln}
Let $\pi\in(\Omega'_N)_0$. The representation
$\mathcal{HP}^{G'_N}_{\overline{G'}_N}\pi^{\oplus}$ 
is a subrepresentation of $\mathcal{P}^{G'_N}_{\overline{G'}_N}\pi^\oplus$ with multiplicity $1$.
\begin{proof}

By \cite[Lemma 1.11]{BKSln} , for any $\tilde{\pi}\in (\Omega_N)_0$  it holds
\[\tilde{\pi}^{K^+_N}=\tilde{\pi}^{{K^+_N}'}.\]
Hence 
\begin{align}\label{parRescomm}
    \mathcal{P}^{G'_N}_{\overline{G'}_N}\pi^\oplus=&
    Res^{G'_N}_{K'_N}(\pi^\oplus)^{{K'_N}^+} &\text{by \eqref{plusrestrictiongln}}\nonumber\\
    &=Res^{G'_N}_{K'_N}(Res^{G_N}_{G'_N}R(\pi))^{{K'_N}^+} &\text{by \cite[Lemma 1.11]{BKSln}} \nonumber\\
    &=Res^{G_N}_{K'_N}(R(\pi))^{K^+_N}& \nonumber\\
    &=Res^{K_N}_{K'_N}(Res^{G_N}_{K_N}(R(\pi))^{K^+_N})=Res^{\overline{G}_N}_{\overline{G'}_N}\mathcal{P}^{G_N}_{\overline{G}_N}(R(\pi)).&
\end{align}

The representation $\mathcal{HP}^{G_N}_{\overline{G}_N}(R(\pi))$ is an irreducible constituent of multiplicity $1$ of $\mathcal{P}^{G_N}_{\overline{G}_N}(R(\pi))$ by Theorem \ref{prephead}.  \\
Therefore $\mathcal{HP}^{G'_N}_{\overline{G'}_N}\pi^{\oplus}=Res^{\overline{G}_N}_{\overline{G'}_N}\mathcal{HP}^{G_N}_{\overline{G}_N}R(\pi)$ is a subrepresentation of  $\mathcal{P}^{G'_N}_{\overline{G'}_N}\pi^\oplus=Res^{\overline{G}_N}_{\overline{G'}_N}\mathcal{P}^{G_N}_{\overline{G}_N}(R(\pi))$.

To prove that it has multiplicity $1$, it is enough to show that there are no irreducible constituents of $\mathcal{P}^{G_N}_{\overline{G}_N}(R(\pi))$ different from $\mathcal{HP}^{G_N}_{\overline{G}_N}R(\pi)$ having its same restriction  to $\overline{G'}_N$. Indeed the restrictions of irreducible representations of $\overline{G}_N$ to $\overline{G'}_N$ are either equal or don't have any common irreducible constituent. 

Irreducible representations of $\overline{G}_N$ have the same restriction to $\overline{G'}_N$ only if they differ by the tensor product with a character  induced by the determinant map. Hence it is sufficient to show that for any $\chi\in \widehat{k_F^*}$ the representation $\mathcal{HP}^{G_N}_{\overline{G}_N}R(\pi)\otimes(\chi\circ det)$ is a subrepresentation of $\mathcal{P}^{G_N}_{\overline{G}_N}(R(\pi))$ if and only if it is  equal to  $\mathcal{HP}^{G_N}_{\overline{G}_N}R(\pi)$.
By Corollary \ref{twistedparrescor}, $\mathcal{HP}^{G_N}_{\overline{G}_N}(R(\pi)\otimes(\chi\circ \det))=\overline{\sigma}_{{\chi}\overline{\Lambda_{R(\pi)}}}$. If $\overline{\sigma}_{{\chi}\overline{\Lambda_{R(\pi)}}}$ is an irreducible constituent of $\mathcal{P}^{G_N}_{\overline{G}_N}(R(\pi))$, then by Theorem \ref{prephead} it holds ${\chi}\overline{\Lambda_{R(\pi)}}(\overline{\tau})\leq \overline{\Lambda_{R(\pi)}}(\overline{\tau})$ for any $\overline{\tau}\in\overline{\mathcal{C}}$, hence $\overline{\Lambda_{R(\pi)}}(\chi^{-1}\otimes\overline{\tau})\leq \overline{\Lambda_{R(\pi)}}(\overline{\tau})$. 
Iterating, we get $\overline{\Lambda_{R(\pi)}}(\chi^{i-1}\otimes\overline{\tau})\leq \overline{\Lambda_{R(\pi)}}(\chi^{i}\otimes\overline{\tau})$ for $i\geq 0$. By transitivity it implies 
\[\overline{\Lambda_{R(\pi)}}(\chi^{-1}\otimes\overline{\tau})\leq \overline{\Lambda_{R(\pi)}}(\overline{\tau})\leq  \overline{\Lambda_{R(\pi)}}(\chi^{q-2}\otimes\overline{\tau})=\overline{\Lambda_{R(\pi)}}(\chi^{-1}\otimes\overline{\tau}),\]
that is ${\chi}\overline{\Lambda_{R(\pi)}}(\overline{\tau})=\overline{\Lambda_{R(\pi)}}(\overline{\tau})$ for any $\overline{\tau}\in\overline{\mathcal{C}}$. Hence if $\mathcal{HP}^{G_N}_{\overline{G}_N}(R(\pi))\otimes(\chi\circ \det)$ is a subrepresentation of $\mathcal{P}^{G_N}_{\overline{G}_N}(R(\pi))$ it holds
\[\mathcal{HP}^{G_N}_{\overline{G}_N}R(\pi)\otimes(\chi\circ \det)=\overline{\sigma}_{{\chi}\overline{\Lambda_{R(\pi)}}}=\overline{\sigma}_{\overline{\Lambda_{R(\pi)}}}=\mathcal{HP}^{G_N}_{\overline{G}_N}R(\pi).\]
\end{proof}
\end{lemma}

\subsubsection{Compatibility of $\mathcal{L}'_N$ and $\mathcal{M}'_N$}

Let
\begin{align}\label{Pi}
    &\Sigma_N:=\{ \pi^\oplus \ |  \ \pi\in \Omega'_N\} &\text{and} \hspace{1cm}\overline{\Sigma}_N:=\{ \overline\pi^\oplus \ |  \ \overline\pi\in \overline{\Omega_N'}\}.
\end{align}
\nomenclature{$\Sigma_N$}{Set of the $\pi^\oplus$ \eqref{Pi}}
\nomenclature{$\overline{\Sigma}_N$}{Set of the $\overline{\pi}^\oplus$ \eqref{Pi}}
By \eqref{plusmaclang} these sets are in bijection with the sets of the fibers of $\mathcal{L}'_N$ and $\mathcal{M}'_N$ respectively.

In other words, we can define the maps
\begin{align}\label{Pimap}
    \mathcal{S}: (\Omega'_N)_0&\rightarrow \Sigma_N &\overline{\mathcal{S}}: \overline{\Omega_N'}\rightarrow \overline\Sigma_N\\
    \pi&\mapsto \pi^\oplus
     &\overline\pi\mapsto \overline\pi^\oplus
\end{align}

\nomenclature{$\mathcal{S}$}{Map $\pi\mapsto\pi^\oplus$ \eqref{Pimap}}
\nomenclature{$\overline{\mathcal{S}}$}{Map $\overline\pi\mapsto\overline\pi^\oplus$ \eqref{Pimap}}

Then the maps $\mathcal{R}$ \eqref{R} and $\mathcal{L}'_N$ factor through $\mathcal{S}$, and we write $\mathcal{R}=(\mathcal{R})_{\oplus}\circ \mathcal{S}$ and $\mathcal{L}'_N= (\mathcal{L}'_N)_{\oplus}\circ \mathcal{S}$, where $(\mathcal{R})_{\oplus}$ and $(\mathcal{L}'_N)_{\oplus}$ are bijections. 
 Analogously, the map $\overline{\mathcal{R}}$ \eqref{Rbar} and $\mathcal{M}'_N$ factor through $\overline{\mathcal{S}}$, and we write  $\overline{\mathcal{R}}=(\overline{\mathcal{R}})_{\oplus}\circ \overline{\mathcal{S}}$ and $\mathcal{M}'_N= (\mathcal{M}'_N)_{\oplus}\circ \overline{\mathcal{S}}$, where $(\overline{\mathcal{R}})_{\oplus}$ and $(\mathcal{M}'_N)_{\oplus}$ are bijections.
 
 With this notation, the following diagrams are commutative:
\begin{equation}
\begin{tikzcd}
&\faktor{(\Omega_N)_0}{\widehat{F^*}} \arrow[r, "\overline{\mathcal{L}_N}"] & \faktor{(\Phi_N)_0}{\widehat{W_F}} \arrow[d, "\eta^{*}"] \\
 & \Sigma_N \arrow[u, "(\mathcal{R})_\oplus"'] \arrow[r, "(\mathcal{L}_N')_\oplus"]    & (\Phi'_N)_0 \\
 (\Omega'_N)_0 \arrow[ru, "\mathcal{S}" description] \arrow[rru, "\mathcal{L}_N'" description] \arrow[ruu, "\mathcal{R}" description] & &   
\end{tikzcd}
\end{equation}

\begin{equation}
\begin{tikzcd}
 &\faktor{(\overline{\Omega_N})}{\widehat{k^*_F}} \arrow[r, "\overline{\mathcal{M}_N}"] &\faktor{\bigg(\faktor{(\Phi_N)_0}{\sim_{I_F}}\bigg)}{\widehat{k_F^*}}\arrow[d, "\eta^{*}"]  \\
 &\overline{\Sigma}_N \arrow[u, "(\overline{\mathcal{R}})_{\oplus}"'] \arrow[r, "(\mathcal{M}_N')_{\oplus}"]  & \faktor{(\Phi'_N)_0}{\sim_{I_F}}     \\
\overline{\Omega'_N} \arrow[ruu, "\overline{\mathcal{R}}" description] \arrow[ru, "\overline{\mathcal{S}}" description] \arrow[rru, "\mathcal{M}_N'" description] &  &
\end{tikzcd}
\end{equation}

\begin{theorem}\label{compatibilitysln}
For any $\pi\in(\Omega'_N)_0$ , it holds
\[(\mathcal{L}'_N)_\oplus\pi^\oplus=(\mathcal{M}'_N)_{\oplus}(\mathcal{HP}^{G'_N}_{\overline{G'}_N}\pi^\oplus)\]
\begin{proof}
By\eqref{plusrestrictiongln}, the restriction from $G_N$ to $G'_N$ (respectively from $\overline{G}_N$ to $\overline{G'}_N$) lands in $\Sigma_N$ (respectively $\overline{\Sigma}_N$), and by \eqref{plusrestrictiongln} it factors through the inverse of the map $(\mathcal{R})_\oplus$ (respectively $(\overline{\mathcal{R}})_{\oplus}$). So the following diagrams are commutative:
\nomenclature{$(f)_{\oplus}$ }{Map such that $f=(f)_\oplus\circ \mathcal{S}$ (or $f=(f)_\oplus\circ \overline{\mathcal{S}}$) }

\begin{equation}\label{diag1}
    \begin{tikzcd}
(\Omega_N)_0 \arrow[r, "\faktor{}{\widehat{F^*}}"] \arrow[d, "\mathcal{L}_N"] \arrow[rr, "Res^{G_N}_{G'_N}", bend left] & \faktor{(\Omega_N)_0}{\widehat{F^*}} \arrow[d, "\overline{\mathcal{L}_N}"] \arrow[r, "Res^{G_N}_{G'_N}"] & \Sigma_N \arrow[d, "(\mathcal{L}'_N)_\oplus"] \\
(\Phi_N)_0 \arrow[r, "\faktor{}{\widehat{W_F}}"] \arrow[rr, "\eta^{*}", bend right]                                       & \faktor{(\Phi_N)_0}{\widehat{W_F}} \arrow[r, "\eta^{*}"]                                                   & (\Phi'_N)_0                  
\end{tikzcd}
\end{equation}

\begin{equation}\label{diag2}
    \begin{tikzcd}
\overline{\Omega_N} \arrow[r, "\faktor{}{\widehat{k_F^*}}"] \arrow[d, "\mathcal{M}_N"] \arrow[rr, "Res^{\overline{G}_N}_{\overline{G'}_N}", bend left] & \faktor{\overline{\Omega_N}}{\widehat{k_F}^*} \arrow[d, "\overline{\mathcal{M}_N}"] \arrow[r, "Res^{\overline{G}_N}_{\overline{G'}_N}"] & \overline{\Sigma}_N \arrow[d, "(\mathcal{M}'_N)_\oplus"] \\
\faktor{(\Phi_N)_0}{\sim_{I_F}} \arrow[r, "\faktor{}{\widehat{k_F^*}}"] \arrow[rr, "\eta^{*}", bend right] & \faktor{\big(\faktor{(\Phi_N)_0}{\widehat{k_F^*}}\big)}{\sim_{I_F}} \arrow[r, "\eta^{*}"]& (\Phi'_N)_0 
\end{tikzcd}
\end{equation}
and so for a representation $\pi\in (\Omega_N')_0$, and for any $R(\pi)\in(\Omega_N)_0$ containing $\pi$ as $G'_N$-subrepresentation, we have
\begin{align*}
(\mathcal{M}'_N)_{\oplus}(\mathcal{HP}^{G'_N}_{\overline{G'}_N}\pi^{\oplus})=&
(\mathcal{M}'_N)_{\oplus} (Res^{\overline{G}_N}_{\overline{G'}_N}\mathcal{HP}^{G_N}_{\overline{G}_N}R(\pi)) &\text{by \eqref{diag2}}\\
&=\eta^{*} \circ \mathcal{M}_N(\mathcal{HP}^{G_N}_{\overline{G}_N}R(\pi))&\text{by Theorem \ref{compatibility}} \\ 
&=\eta^{*} \circ \mathcal{L}_N(R(\pi)) &\text{by \eqref{diag1}}\\
&=(\mathcal{L}'_N)_\oplus (Res^{G_N}_{G'_N} R(\pi))= (\mathcal{L}'_N)_\oplus \pi^\oplus. &\text{by \eqref{plusrestrictiongln}}\\
\end{align*}
\end{proof}
    
\end{theorem}

\begin{rmk}
Theorem \ref{compatibilitysln} does not depend on the choice of the maximal compact subgroup $K'_N$ of $G'_N$.\\
Indeed, any maximal compact subgroup $J'_N$ of $G'_N$ is of the form $J'_N={}^{g}K'_N$ for some $g\in G_N$.
For any $\pi\in(\Omega'_N)_0$, we have ${}^g\pi^{\oplus}=\pi^{\oplus}$ and so
\begin{align*}
    (Res^{G'_N}_{J'_N}\pi^\oplus)^{{J'_N}^+}=
    (Res^{G'_N}_{{}^gK'_N}\pi^\oplus)^{{}^g{K'_N}^+}=
({}^g Res^{G'_N}_{K'_N} ({}^{g^{-1}}\pi^\oplus))^{{}^g{K'_N}^+}= 
{}^{g}((Res^{G'_N}_{K'_N}\pi^\oplus)^{{K'_N}^+}).
\end{align*}
 We have  $\faktor{J_N'}{{J_N'}^+}\cong \faktor{K_N'}{{K_N'}^+}\cong{\overline{G'}_N}$ where the first isomorphism is induced by conjugation by $g$. Therefore $(Res^{G'_N}_{J'_N}\pi^\oplus)^{{J'_N}^+}$ and $(Res^{G'_N}_{K'_N}\pi^\oplus)^{{K'_N}^+}$ are isomorphic as $\overline{G'}_N$ representations.
\end{rmk}

\subsection{Enhancement of the compatibility via the parametrization of the fibers}
The central result of this section is Theorem \ref{finalcomp}, where Theorem \ref{compatibilitysln} is refined in order to consider representations individually, rather than fibers.

For the rest of this section, we fix $(\rho,E)\in  (\Phi'_N)_0$, and we let $\pi_{(\rho,E)}\in\mathcal{L'}^{-1}(\rho,E)$ and  $\overline{\pi}_{(\rho,E)}\in\mathcal{M'}^{-1}((\rho,E)_{I_F})$. We fix $(\tilde{\rho},E)\in  (\Phi_N)_0$  such that $\eta\circ\tilde{\rho}=\rho$, and we set $\tilde{\pi}_{(\tilde{\rho},E)}:=\mathcal{L}_N^{-1}(\tilde{\rho},E)$ and $\widetilde{\overline{\pi}}_{(\tilde{\rho},E)}:=\mathcal{M}_N^{-1}((\tilde{\rho},E)_{I_F})$. 

By \cite[Theorem 4.3]{Gel}, the character group  $A_{PGL_N(\mathbb{C})}(\rho,E)^{\wedge}$ acts simply transitively on  
${\mathcal{L}'_N}^{-1}(\rho,E)$, and by Theorem \ref{fiber2} the character group  $A_{PGL_N(\mathbb{C})}(\rho|_{I_F},\rho(Fr)C^0_{PGL_N(\mathbb{C})}(\rho|_{I_F}),E)^{\wedge}$ acts simply transitively on  
${\mathcal{M}'_N}^{-1}(\rho,E)$.
We now describe explicitly these actions.

By \cite[Theorem 4.3]{Gel}, there is an isomorphism 
\begin{align}\label{xif}
  \Xi_{F}: A_{PGL_N(\mathbb{C})}(\rho,E)&\rightarrow stab_{\widehat{W_F}}(\tilde{\rho},E)\\
  [g]&\mapsto \chi_g:=\tilde{g}\tilde{\rho}\tilde{g}^{-1}\tilde{\rho}^{-1},\nonumber
\end{align} 
\nomenclature{$\Xi_F$}{Isomorphism \eqref{xif}}
where $\tilde{g}$ is a lift of $g$ in $GL_N(\mathbb{C})$. 
\begin{rmk}\label{tameenough}
In \cite[Theorem 4.3]{Gel} $\Xi_F$ is defined from $A_{PGL_N(\mathbb{C})}(\rho,E)$ to the stabilizer in the whole character group of $W_F$. But a character $\chi$ of $W_F$ stabilizing a tame Langlands parameter $(\rho, E)\in (\Phi_N)_0$ is necessarily tame: indeed since $(\rho,E)$ is tame, $\rho|_{P_F}=1$, hence
\[\chi|_{P_F}=\chi|_{P_F}\otimes\rho|_{P_F}=(\chi\otimes\rho)|_{P_F}\cong \rho|_{P_F}=1 .\]
\end{rmk}

The bijection of class field theory induces an identification $stab_{\widehat{W_F}}(\tilde{\rho},E)\cong stab_{\widehat{F^{*}}}(\tilde{\pi}_{(\tilde{\rho},E)})$. We still denote by $\Xi_F$ the isomorphism $A_{PGL_N(\mathbb{C})}(\rho,E)\cong stab_{\widehat{F^{*}}}(\tilde{\pi}_{(\tilde{\rho},E)})$ obtained by this identification.
Hence $A_{PGL_N(\mathbb{C})}(\rho,E)$
is a finite abelian group and $\Xi_{F}$ induces the dual  isomorphism  
\begin{align*}
  \widehat{\Xi}_{F}: stab_{\widehat{F^*}}( \tilde{\pi}_{(\tilde{\rho},E)})^{\wedge} \rightarrow
  A_{PGL_N(\mathbb{C})}(\rho,E)^{\wedge}.
\end{align*}
By \cite[Corollary 2.2]{Gel} it holds $stab_{\widehat{F^*}}( \tilde{\pi}_{(\tilde{\rho},E)})\cong\big(\faktor{F^*}{stab_{F^*}(\pi_{(\rho, E)})}\big)^{\wedge}$. Dualizing we have a canonical isomorphism 
\begin{align}\label{Psi}
   \Psi: \faktor{F^*}{stab_{F^*}(\pi_{(\rho, E)})}&\rightarrow stab_{\widehat{F^*}}( \tilde{\pi}_{(\tilde{\rho},E)})^{\wedge}\\
    x&\mapsto \psi_{x}\nonumber
\end{align}
\nomenclature{$\Psi$}{Isomorphism \eqref{Psi}}
defined by $\psi_x(\chi):=\chi(x)$ for any $\chi\in stab_{\widehat{F^*}}( \tilde{\pi}_{(\tilde{\rho},E)})$. Then $ \widehat{\Xi}_{F}\circ \Psi$ is a canonical isomorphism between  $\faktor{F^*}{stab_{F^*}(\pi_{(\rho, E)})}$ and $A_{PGL_N(\mathbb{C})}(\rho,E)^{\wedge}$. The determinant induces an isomorphism $det:\faktor{G_N}{stab_{G_N}(\pi_{(\rho,E)})}\rightarrow \faktor{F^*}{stab_{F^*}(\pi_{(\rho, E)})}$. By abuse of notation we denote with the same symbol elements corresponding to each other through this isomorphism. For any $\psi\in A_{PGL_N(\mathbb{C})}(\rho,E)^{\wedge}$, let $x_{\psi}=(\widehat{\Xi}_{F}\circ \Psi)^{-1}(\psi)$. The action of  $A_{PGL_N(\mathbb{C})}(\rho,E)^{\wedge}$  on ${\mathcal{L}'}^{-1}(\rho, E)$ is given by $\psi.\pi={}^{x_{\psi}}\pi$. 

The canonical isomorphism $\Xi$ defined in Theorem \ref{fiber2} induces by duality the isomorphism  
\begin{align*}
 \widehat{\Xi}: stab_{\widehat{k_{F^*}}}(\tilde{\rho},E)^{\wedge} \rightarrow   A_{PGL_N(\mathbb{C})}(\rho|_{I_F},\rho(Fr)C^0_{PGL_N(\mathbb{C})}(\rho|_{I_F}),E)^{\wedge}
\end{align*}
and by Lemma \ref{fiber0} there is a canonical isomorphism
\begin{align}\label{Psibar}
   \overline{\Psi}: \faktor{k_F^*}{stab_{k_F^*}(\overline{\pi}_{(\rho,E)})}&\rightarrow stab_{\widehat{k_F^*}}( \widetilde{\overline{\pi}}_{(\tilde{\rho},E)})^{\wedge}\\
    x&\mapsto \overline{\psi}_{x}\nonumber
\end{align}
\nomenclature{$\overline{\Psi}$}{Isomorphism \eqref{Psibar}}
defined by $\overline{\psi}_x(\overline\chi):=\overline\chi(x)$ for any $\overline{\chi}\in stab_{\widehat{k_F^*}}( \tilde{\pi}_{(\tilde{\rho},E)})$. As before, the determinant induces an isomorphism $\faktor{\overline{G}_N}{stab_{\overline{G}_N}(\overline{\pi}_{(\rho,E)})}\cong \faktor{k_F^*}{stab_{k_F^*}(\overline{\pi}_{(\rho,E)})}$, and we identify corresponding elements in these groups. For any $\overline{\psi}\in   A_{PGL_N(\mathbb{C})}(\rho|_{I_F},\rho(Fr)C^0_{PGL_N(\mathbb{C})}(\rho|_{I_F}),E)$, let $x_{\overline{\psi}}=(\widehat{\Xi}\circ \overline{\Psi})^{-1}(\overline{\psi})$. The action of $A_{PGL_N(\mathbb{C})}(\rho|_{I_F},\rho(Fr)C^0_{PGL_N(\mathbb{C})}(\rho|_{I_F}),E)^{\wedge}$ is given by $\overline{\psi}.\pi={}^{x_{\overline{\psi}}}\pi$.

Recall that the group morphism \eqref{inclusioncomp} induces by duality a canonical group morphism
\[\hat{\iota}:A_{PGL_N(\mathbb{C})}(\rho|_{I_F},\rho(Fr)C^0_{PGL_N(\mathbb{C})}(\rho|_{I_F}),E)^{\wedge}\rightarrow A_{PGL_N(\mathbb{C})}(\rho, E)^{\wedge}.\]

\begin{lemma}\label{comm0}
 The following diagram commutes:
  
 \begin{equation}\label{bottom}
\begin{tikzcd}
stab_{\widehat{k_F^*}}(\widetilde{\overline{\pi}}_{\tilde{\rho};E)}) \arrow[r, "(Res|_{\of{F}^*})^{\wedge}"] \arrow[d, "\widehat{\Xi}"]   & stab_{\widehat{F^*}}({\tilde{\pi}}_{\tilde{\rho};E)})^{\wedge} \arrow[d, "\widehat{\Xi}_{F}"] \\
{ A_{PGL_N(\mathbb{C})^{\wedge}}(\rho|_{I_F},\rho(Fr)C^0_{PGL_N(\mathbb{C})}(\rho|_{I_F}),E)^{\wedge}} \arrow[r, "\hat{\iota}"]                & {A_{PGL_N(\mathbb{C})}(\rho, E)^{\wedge}}                                                 
\end{tikzcd}
 \end{equation}

\begin{proof}
The diagram below,  where the vertical maps labelled by $\cong$ are the identifications given by local class field theory, commutes by definition of the maps $\Xi$ \eqref{chibar}
 and $\Xi_F$ \eqref{xif}:
     \begin{equation*}
\begin{tikzcd}
stab_{\widehat{k_F^*}}(\widetilde{\overline{\pi}}_{\tilde{\rho};E)})                                           &  & stab_{\widehat{F^*}}({\tilde{\pi}}_{\tilde{\rho};E)}) \arrow[ll, "Res|_{\of{F}^*}"] \\
{stab_{\widehat{k_F^*}}(\tilde{\rho},E)} \arrow[u, "\cong"]                                                    &  & stab_{\widehat{W_F}}(\tilde{\rho};E) \arrow[ll, "Res|_{I_F}"] \arrow[u, "\cong"]     \\
{ A_{PGL_N(\mathbb{C})}(\rho|_{I_F},\rho(Fr)C^0_{PGL_N(\mathbb{C})}(\rho|_{I_F}),E)^{\wedge}} \arrow[u, "\Xi"] &  & {A_{PGL_N(\mathbb{C})}(\rho, E)} \arrow[ll, "\iota"] \arrow[u, "\Xi_F"]             
\end{tikzcd}
 \end{equation*} 
 Dualizing, we obtain the statement.
\end{proof}

\end{lemma}

Any section of the projection $\of{F}^*\rightarrow \faktor{\of{F}^*}{1+\pf{F}}\cong k_F$ induces a map \[\widetilde{\mathcal{J}}:k_F^*\rightarrow \faktor{F^*}{stab_{F^*}(\pi_{(\rho,E)})}\]
This map does not depend on the choice of the section: in order to check that $(1+\pf{F})\leq stab_{F^*}(\pi_{(\rho,E)})$ it is enough to check that $\Psi(1+\pf{F})=1$, because $\Psi$ is a group isomorphism. For any $y\in (1+\pf{F})$ it holds  ${\Psi}(y)(\chi)=\chi(y)=1$ since any $\chi\in  stab_{\widehat{F^*}}( \tilde{\pi}_{(\tilde{\rho},E)})$ is tame (see Remark \ref{tameenough}). 

\begin{lemma}\label{comm1}
With notation as above, the map $\widetilde{\mathcal{J}}$ induces a map
\begin{align}\label{J}
\mathcal{J}:\faktor{k_F^*}{stab_{k_F^*}(\overline{\pi}_{(\rho,E)})}&\rightarrow \faktor{F^*}{stab_{F^*}(\pi_{(\rho,E)})}
\end{align}
\nomenclature{$\mathcal{J}$}{The map in \eqref{J}}
that makes the following diagram commutative:
\begin{equation}\label{top}
\begin{tikzcd}
\faktor{k_F^*}{stab_{k_F^*}(\overline{\pi}_{(\rho, E)})} \arrow[r, "\mathcal{J}"] \arrow[d, "\overline{\Psi}"]  & \faktor{F^*}{stab_{F^*}(\pi_{(\rho,E)})} \arrow[d, "\Psi"] \\
stab_{\widehat{k_F^*}}(\widetilde{\overline{\pi}}_{\tilde{\rho};E)})^{\wedge} \arrow[r, "(Res|_{\of{F}^*})^{\wedge}"] &  stab_{\widehat{F^*}}({\tilde{\pi}}_{\tilde{\rho};E)})^{\wedge}      
\end{tikzcd}
 \end{equation}
\begin{proof}

Let $p: {k_F^*}\rightarrow \faktor{k_F^*}{stab_{k_F^*}(\overline{\pi}_{(\rho,E)})}$ denote the natural projection. It holds
\begin{equation}\label{coltilde}
    (Res|_{{\of{F}}^*})^{\wedge}\circ \overline{\Psi}\circ p=\Psi\circ\widetilde{\mathcal{J}}.
\end{equation}
Indeed for any  $x\in {k_F^*}$ and for any $\chi\in\widehat{F^*}$
\begin{align*}
    (Res|_{\of{F}^*})^{\wedge}\circ \overline{\Psi}\circ p (x)(\chi)&= (Res|_{\of{F}^*})^{\wedge}(\overline{\psi}_{p(x)})(\chi)=\overline{\psi}_{p(x)}(\chi|_{\of{F}^*})\\&=\chi|_{\of{F}^*}(x)=\chi(\widetilde{\mathcal{J}}(x))=\psi_{\widetilde{\mathcal{J}}(x)}(\chi)=\Psi\circ\widetilde{\mathcal{J}}(x)(\chi).
\end{align*}

Therefore ${stab_{k_F^*}(\overline{\pi}_{(\rho,E)})}\leq Ker(\widetilde{\mathcal{J}})$ because $\Psi$ is an isomorphism and ${stab_{k_F^*}(\overline{\pi}_{(\rho,E)})}\leq Ker(p)$. 
 So $\widetilde{\mathcal{J}}$ induces a map $\mathcal{J}$ such that $\widetilde{\mathcal{J}}={\mathcal{J}}\circ p$. Since $p$ is an epimorphism,  equation \eqref{coltilde} yields
 \[ (Res|_{{\of{F}}^*})^{\wedge}\circ \overline{\Psi}=\Psi\circ\mathcal{J}\]
 that is the commutativity of the diagram in the statement. 
\end{proof}
\end{lemma}

\begin{prop}\label{comm}
   With notation as above, the following diagram commutes:
   \begin{equation}
\begin{tikzcd}
\faktor{k_F^*}{stab_{k_F^*}(\overline{\pi}_{(\rho, E)})} \arrow[r, "\mathcal{J}"] \arrow[d, "\widehat{\Xi}\circ\overline{\Psi}"]   & \faktor{F^*}{stab_{F^*}(\pi_{(\rho,E)})} \arrow[d, "\widehat{\Xi}_F\circ\Psi"] \\
{ A_{PGL_N(\mathbb{C})}(\rho|_{I_F},\rho(Fr)C^0_{PGL_N(\mathbb{C})}(\rho|_{I_F}),E)^{\wedge}} \arrow[r, "\hat{\iota}"]                & {A_{PGL_N(\mathbb{C})}(\rho, E)^{\wedge}}                                            
\end{tikzcd}
   \end{equation}
\begin{proof}
The statement is obtained stacking the commutative diagram \eqref{top} in Lemma \ref{comm0} on top of the commutative diagram \eqref{bottom} in Lemma \ref{comm1}. 
 \end{proof}
\end{prop}

\begin{prop}\label{equiP}
 Let $\overline{\psi}\in A_{PGL_N(\mathbb{C})}(\rho|_{I_F},\rho(Fr)C^0_{PGL_N(\mathbb{C})}(\rho|_{I_F}),E)^{\wedge}$. If $\overline{\pi}_{(\rho,E)}$ is an irreducible constituent of $\mathcal{P}^{G'_N}_{\overline{G'}_N}\pi_{(\rho,E)}$, then $\overline{\psi}.\overline{\pi}_{(\rho,E)}$ is an irreducible constituent of $\mathcal{P}^{G'_N}_{\overline{G'}_N}(\hat{\iota}(\overline{\psi}).\pi_{(\rho,E)})$.
\begin{proof}
Let $x=(\widehat{\Xi}\circ \overline{\Psi})^{-1}(\overline{\psi})$ so that $\overline{\psi}.\overline{\pi}_{(\rho,E)}={}^x\overline\pi_{(\rho,E)}$. By Proposition \eqref{comm}
\[(\widehat{\Xi}_{F}\circ \Psi)^{-1}(\hat{\iota}(\overline{\psi}))= \mathcal{J}((\widehat{\Xi}\circ \overline{\Psi})^{-1}(\overline{\psi}))= \mathcal{J}(x)\]
and therefore $\hat{\iota}(\overline{\psi}).\pi_{(\rho,E)}={}^{\mathcal{J}(x)}\pi_{(\rho,E)}$.

We need to prove that if $\overline{\pi}_{(\rho,E)}$ is an irreducible constituent of $\mathcal{P}^{G'_N}_{\overline{G'}_N}\pi_{(\rho,E)}$, then  $\overline{\psi}.\overline{\pi}_{(\rho,E)}={}^x\overline\pi_{(\rho,E)}$ is an irreducible constituent of $\mathcal{P}^{G'_N}_{\overline{G'}_N}(\hat{\iota}(\overline{\psi}).\pi_{(\rho,E)})=\mathcal{P}^{G'_N}_{\overline{G'}_N}({}^{\mathcal{J}(x)}\pi_{(\rho,E)})$. Rephrasing it, we need to show that $Hom_{\overline{G'}_N}(\overline\pi_{(\rho,E)},\mathcal{P}^{G'_N}_{\overline{G'}_N}(\pi_{(\rho,E)}))\neq0$ implies that 
$Hom_{\overline{G'}_N}({}^x\overline\pi_{(\rho,E)},\mathcal{P}^{G'_N}_{\overline{G'}_N}({}^{\mathcal{J}(x)}\pi_{(\rho,E)}))\neq0$. 

The inflation by the compact subgroup ${K'_N}^+$ is left adjoint to taking the ${K'_N}^+$-fixed subspace, so 
\begin{equation}\label{adjinfl}
    Hom_{\overline{G'}_N}({}^x\overline{\pi}_{(\rho,E)},\mathcal{P}^{G'_N}_{\overline{G'}_N}({}^{\mathcal{J}(x)}\pi_{(\rho,E)}))\cong Hom_{K'_N}(Infl_{\overline{G'}_N}^{K'_N}{}^x\overline\pi_{(\rho,E)}, Res^{G'_N}_{K'_N}({}^{\mathcal{J}(x)}\pi_{(\rho,E)})).
\end{equation}
The lifting of the map $\mathcal{J}$ through the map $det$, that by abuse of notation we still denote by  $\mathcal{J}$, is the map $\mathcal{J}: \faktor{\overline{G}_N}{stab_{\overline{G}_N}(\overline{\pi}_{(\rho,E)})}\rightarrow \faktor{{G_N}}{Stab_{G_N}(\pi_{(\rho,E)})}$ induced by any section of the projection map $r:K_N\rightarrow \faktor{K_N}{K_N^+}\cong \overline{G}_N$, so the element $\mathcal{J}(x)\in \faktor{{G_N}}{Stab_{G_N}(\pi)}$ has a representative in $K_N\leq G_N$. In the following of this proof we fix such a representative. The final result does not depend on this choice, because ${}^{\mathcal{J}(x)}\pi_{(\rho,E)}$ does not depend on the choice of the representative of $\mathcal{J}(x)$ in $G_N$.

We observe that
\begin{equation}\label{portafuorix}
    Infl_{\overline{G'}_N}^{K'_N}({}^x\overline\pi_{(\rho,E)})={}^{\mathcal{J}(x)}Infl_{\overline{G'}_N}^{K'_N}\overline{\pi}_{(\rho,E)},
\end{equation}
where the right-hand side is a representation of $K'_N$ by normality of $K'_N$ in $K_N$. Indeed, since the natural projection $K'_N\rightarrow \faktor{K'_N}{{K'_N}^+}$ is obtained by restriction of $r$, for any $k\in K'$ 
\begin{align*}
    Infl_{\overline{G'}_N}^{K'_N}({}^x\overline\pi_{(\rho,E)})(k)=\overline{\pi}_{(\rho,E)})(x^{-1}r(k)x)=\overline{\pi}_{(\rho,E)})(r(\mathcal{J}(x)^{-1}k\mathcal{J}(x)))=
{}^{\mathcal{J}(x)}Infl_{\overline{G'}_N}^{K'_N}\overline{\pi}_{(\rho,E)}(k).
\end{align*}

Combining equation \eqref{adjinfl} and \eqref{portafuorix} gives
\begin{align*}
    Hom_{\overline{G'}_N}({}^x\overline{\pi}_{(\rho,E)},\mathcal{P}^{G'_N}_{\overline{G'}_N}({}^{\mathcal{J}(x)}\pi_{(\rho,E)}))&\cong
     Hom_{K'_N}({}^{\mathcal{J}(x)}Infl_{\overline{G'}_N}^{K'_N}\overline\pi_{(\rho,E)},Res^{G'_N}_{K'_N}({}^{\mathcal{J}(x)}\pi_{(\rho,E)}))\\
   &= Hom_{K'_N}(Infl_{\overline{G'}_N}^{K'_N}\overline\pi_{(\rho,E)}, {}^{{\mathcal{J}(x)}^{-1}}Res^{G'_N}_{K'_N}({}^{\mathcal{J}(x)}\pi_{(\rho,E)})).
\end{align*}
We rewrite ${}^{\mathcal{J}(x)^{-1}}Res^{G'_N}_{K'_N}({}^{\mathcal{J}(x)}\pi_{(\rho,E)})=Res^{G'_N}_{{}^{\mathcal{J}(x)^{-1}}K'_N}(\pi_{(\rho,E)})$. By normality of $K'_N$ in $K_N$, it holds $Res^{G'_N}_{{}^{\mathcal{J}(x)^{-1}}K'_N}(\pi_{(\rho,E)})=Res^{G'_N}_{K'_N}(\pi_{(\rho,E)})$. Then
\begin{align*}
    Hom_{\overline{G'}_N}({}^x\overline{\pi}_{(\rho,E)},\mathcal{P}^{G'_N}_{\overline{G'}_N}({}^{\mathcal{J}(x)}\pi_{(\rho,E)}))&\cong Hom_{K'_N}(Infl_{\overline{G'}_N}^{K'_N}\overline\pi_{(\rho,E)}, Res^{G'_N}_{K'_N}\pi_{(\rho,E)})\\
    &\cong Hom_{\overline{G'}_N}(\overline\pi_{(\rho,E)},\mathcal{P}^{G'_N}_{\overline{G'}_N}(\pi_{(\rho,E)}))\neq0.
\end{align*}

\end{proof}
\end{prop}

\begin{lemma}\label{basepointcomp}
For any $\overline{\pmb{\pi}}_{(\rho,E)}\in\mathcal{M'}^{-1}((\rho,E)_{I_F})$ there exists a unique $\pmb{\pi}_{(\rho,E)}\in\mathcal{L'}^{-1}(\rho,E)$ such that $\overline{\pmb{\pi}}_{(\rho,E)}$ is an irreducible constituent of $\mathcal{P}^{G'_N}_{\overline{G'}_N}\pmb{\pi}_{(\rho,E)}$. Moreover $\overline{\pmb{\pi}}_{(\rho,E)}$ has multiplicity $1$ in $\mathcal{P}^{G'_N}_{\overline{G'}_N}\pmb{\pi}_{(\rho,E)}$. 
\begin{proof}
Denote ${\pi}_{(\rho,E)}^{\oplus}=\bigoplus_{\pi\in\mathcal{L'}^{-1}(\rho,E) }{\pi}$ and $\overline\pi_{(\rho,E)}^{\oplus}=\bigoplus_{\overline{\pi}\in\mathcal{M'}^{-1}(\rho,E) }\overline{\pi}$. By Theorem \ref{compatibilitysln}, it holds 
$\mathcal{HP}^{G'_N}_{\overline{G'}_N} \pi^{\oplus}_{(\rho,E)}=\overline\pi_{(\rho,E)}^{\oplus}$. Hence by Lemma \ref{HPsln} $\overline\pi_{(\rho,E)}^{\oplus}$ is a subrepresentation of $\mathcal{P}^{G'_N}_{\overline{G'}_N}(\pi_{(\rho,E)}^{\oplus})$ of multiplicity $1$:  
\[Hom_{\overline{G'}_N}(\overline\pi_{(\rho,E)}^{\oplus},\mathcal{P}^{G'_N}_{\overline{G'}_N}\pi_{(\rho,E)}^{\oplus})=Hom_{\overline{G'}_N}(\overline\pi_{(\rho,E)}^{\oplus}, \overline\pi_{(\rho,E)}^{\oplus}).\]
Since $\overline{\pmb{\pi}}_{(\rho,E)}$ is an irreducible constituent of multiplicity $1$ of $\overline\pi_{(\rho,E)}^{\oplus}$, it is an irreducible constituent of multiplicity $1$ of $\mathcal{P}^{G'_N}_{\overline{G'}_N}(\pi_{(\rho,E)}^{\oplus})=\bigoplus_{{\pi}\in\mathcal{L'}^{-1}(\rho,E)}( \mathcal{P}^{G'_N}_{\overline{G'}_N}\pi)$. Therefore it is an irreducible constituent of multiplicity $1$ of $\mathcal{P}^{G'_N}_{\overline{G'}_N}\pmb{\pi}_{(\rho, E)}$ for a uniquely determined $\pmb{\pi}_{(\rho, E)}\in\mathcal{L'}^{-1}(\rho,E)$. 
\end{proof}
\end{lemma}

From now on, we fix a representation $\overline{\pmb{\pi}}_{(\rho,E)}\in\mathcal{M'}^{-1}((\rho,E)_{I_F})$ and we denote by $\pmb{\pi}_{(\rho,E)}$ the representation in $\mathcal{L'}^{-1}(\rho,E)$ as in Lemma \ref{basepointcomp}.

For any $\psi \in A_{PGL_N(\mathbb{C})}(\rho, E)^{\wedge} $, the preimage ${\hat{\iota}}^{-1}(\psi)$ is either empty or a coset of $Ker(\hat{\iota})$ in $A_{PGL_N(\mathbb{C})}(\rho|_{I_F},\rho(Fr)C^0_{PGL_N(\mathbb{C})}(\rho|_{I_F}),E)^{\wedge}$. We define the head of parahoric restriction for the representations in $\mathcal{L'}^{-1}(\rho,E)$.
\begin{defi}\label{HPSlnSingle}
With notation as above, we define
\[\mathcal{HP}^{G'_N}_{\overline{G'}_N}(\psi.\pmb{\pi}_{(\rho,E)}):=\bigoplus_{\overline{\psi}\in{\hat{\iota}}^{-1}(\psi)}\overline{\psi}.\overline{\pmb{\pi}}_{(\rho,E)}.\]
\end{defi}
\nomenclature{$\mathcal{HP}^{G'_n}_{\overline{G'}_n}$}{Head of parahoric restriction for $G'_n$}
Definition \ref{HPSlnSingle} is compatible with the Definition \ref{HPSLn} because
\[\mathcal{HP}^{G'_N}_{\overline{G'}_N}(\pi_{(\rho,E)}^{\oplus})=\bigoplus _{\psi\in A_{PGL_N(\mathbb{C})}(\rho, E)^{\wedge}}\mathcal{HP}^{G'_N}_{\overline{G'}_N}(\psi.\pmb{\pi}_{(\rho,E)}).\]
\begin{rmk}
Let $\psi \in A_{PGL_N(\mathbb{C})}(\rho, E)^{\wedge} $. By the definition of $\hat{\iota}$ and Frobenius reciprocity, the coset ${\hat{\iota}}^{-1}(\psi)$ in $A_{PGL_N(\mathbb{C})}(\rho|_{I_F},\rho(Fr)C^0_{PGL_N(\mathbb{C})}(\rho|_{I_F}),E)^{\wedge}$ can be described as follows. If $\psi\in A_{PGL_N(\mathbb{C})}(\rho, E)^{\wedge}$  does not factor through $\iota$, then ${\hat{\iota}}^{-1}(\psi)$ is empty. Otherwise, there exists a character $\xi$  of $\iota(A_{PGL_N(\mathbb{C})}(\rho, E))$ such that $\xi\circ\iota=\psi$. Then ${\hat{\iota}}^{-1}(\psi)$ is the set of the irreducible constituents of $ind_{\iota(A_{PGL_N(\mathbb{C})}(\rho, E))}^{A_{PGL_N(\mathbb{C})}(\rho|_{I_F},\rho(Fr)C^0_{PGL_N(\mathbb{C})}(\rho|_{I_F}),E)}\xi$. 
\end{rmk}
\begin{theorem}\label{finalcomp}
  For any $\psi \in A_{PGL_N(\mathbb{C})}(\rho, E)^{\wedge} $, the representation $\mathcal{HP}^{G'_N}_{\overline{G'}_N}(\psi.\pmb{\pi}_{(\rho,E)})$ is a subrepresentation of $\mathcal{P}^{G'_N}_{\overline{G'}_N}(\psi.\pmb{\pi}_{(\rho,E)})$ of multiplicity $1$.
  
    Moreover if $\overline{\pi}\in{\mathcal{M}'}^{-1}((\rho, E)_{I_F})$ is an irreducible component of $\mathcal{P}^{G'_N}_{\overline{G'}_N}(\psi.\pmb{\pi}_{(\rho,E)})$, then $\overline{\pi}$ is an irreducible component of $\mathcal{HP}^{G'_N}_{\overline{G'}_N}(\psi.\pmb{\pi}_{(\rho,E)})$.
\begin{proof}
By Proposition \ref{equiP}, for any $\overline{\psi}\in\hat{\iota}^{-1}(\psi)$ the representation $\overline{\psi}.\overline{\pmb{\pi}}_{(\rho,E)}$ is an irreducible constituent of $\mathcal{P}^{G_N}_{G'_N}\psi.\pmb{\pi}_{(\rho,E)}$, and by Lemma \ref{basepointcomp} it has multiplicity $1$. Hence $\mathcal{HP}^{G'_N}_{\overline{G'}_N}(\psi.\pmb{\pi}_{(\rho,E)})=\bigoplus_{\overline{\psi}\in{\hat{\iota}}^{-1}(\psi)}\overline{\psi}.\overline{\pmb{\pi}}_{(\rho,E)}$ is a subrepresentation of $\mathcal{P}^{G_N}_{G'_N}\psi.\pmb{\pi}_{(\rho,E)}$ of multiplicity $1$. 

If $\overline{\pi}\in{\mathcal{M}'}^{-1}((\rho, E)_{I_F})$, there exists $\xi\in  A_{PGL_N(\mathbb{C})}(\rho|_{I_F},\rho(Fr)C^0_{PGL_N(\mathbb{C})}(\rho|_{I_F}),E)$ such that $\overline{\pi}=\xi.\overline{\pmb{\pi}}_{(\rho,E)}$ by transitivity of the action. By Proposition \ref{equiP} $\overline{\pi}$ is an irreducible constituent of $\mathcal{P}^{G_N}_{G'_N}\hat{\iota}(\xi).\pmb{\pi}_{(\rho,E)}$. Then if $\overline\pi\leq \mathcal{P}^{G_N}_{G'_N}\psi.\pmb{\pi}_{(\rho,E)}$, by the uniqueness statement in Lemma \ref{basepointcomp} it holds $\psi.\pmb{\pi}_{(\rho,E)}=\hat{\iota}(\xi).\pmb{\pi}_{(\rho,E)}$ and since the action of  $A_{PGL_N(\mathbb{C})}(\rho, E)^{\wedge}$ over ${\mathcal{L}'_N}^{-1}(\rho,E)$ is free, $\hat{\iota}(\xi)=\psi$. Therefore, $\xi\in\hat{\iota}^{-1}(\psi)$, and so $\overline\pi=\xi.\overline{\pmb{\pi}}_{(\rho,E)}$ is an irreducible constituent of $\mathcal{HP}^{G'_N}_{\overline{G'}_N}(\psi.\pi_{(\rho,E)})$.
\end{proof}
\end{theorem}

\begin{rmk}
    Theorem \ref{finalcomp} depends on the choice of the maximal compact subgroup $K'_N$ of $G'_N$. When a different maximal compact subgroup $J'_N$ is considered, the association in Lemma \ref{basepointcomp} between a representation $\overline{\pi}_{(\rho,E)}\in{\mathcal{M}'_N}^{-1}(\rho, E)$ and $\pi_{(\rho, E)}\in{\mathcal{L}'_N}^{-1}(\rho, E)$ given by $Hom_{\overline{G'}_N}(\overline{\pi}_{(\rho,E)}, \mathcal{P}^{G'_N}_{\overline{G'}_N}\pi_{(\rho, E)})\neq 0$, for $(\rho,E)\in(\Phi'_N)_0$, varies as follows.
    Let $g\in G_N$ be such that $J_N={}^g K'_N$. Then ${}^g{\pi}_{(\rho, E)}\in{\mathcal{L}'_N}^{-1}(\rho, E)$, and
\[(Res^{G_N}_{J'_N}({}^g\pi_{(\rho,E)}))^{{J'_N}^+}=(Res^{G_N}_{{}^g K'_N}{}({}^g{\pi}_{(\rho,E))})^{{}^g{K'_N}^+}={}^g((Res^{G_N}_{K'_N}\pi_{(\rho,E)})^{{K'_N}^+}).\]
We have $\overline{G'}_N\cong \faktor{K'_N}{{K'_N}^+}\cong \faktor{J'_N}{{J'_N}^+}$ where the last isomorphism is induced by conjugation by $g$. Hence $Res^{G_N}_{J'_N}{}^g\pi_{(\rho,E)}^{{J'_N}^+}$ and $Res^{G_N}_{K'_N}\pi_{(\rho,E)}^{{K'_N}^+}$ are isomorphic as $\overline{G'}_N$ representations. It follows that $Hom_{\overline{G'}_N}(\overline{\pi}_{(\rho,E)}, (Res^{G_N}_{J'_N}{}^g\pi_{(\rho,E)})^{{J'_N}^+})\neq 0$. Therefore, considering the parahoric restriction with respect to $J'_N$ rather than $K'_N$,  Lemma \ref{basepointcomp} would associate to $\overline{\pi}_{(\rho,E)}\in{\mathcal{M}'_N}^{-1}(\rho, E)$ the representation ${}^{g}\pi_{(\rho, E)}\in{\mathcal{L}'_N}^{-1}(\rho, E)$ rather than $\pi_{(\rho, E)}\in{\mathcal{L}'_N}^{-1}(\rho, E)$. 

However, the equivariance result in Proposition \eqref{equiP} still holds when considering parahoric restriction with respect to a different maximal parahoric subgroup. 

We conclude that replacing $K'_N$ by $J'_N={}^g K'_N$ Theorem \ref{finalcomp} holds replacing $\pmb{\pi}_{(\rho, E)}$ with ${}^{g}\pmb{\pi}_{(\rho, E)}$.

\end{rmk}

\section{Acknowledgements}
I would like to express my deepest gratitude to my PhD supervisors, Giovanna Carnovale and  Dan Ciubotaru, for their great guidance and insightful feedback throughout the preparation of this work. Their expertise and mentorship were fundamental to the development and completion of this paper.

I am also sincerely grateful to Anne-Marie Aubert for the enlightening discussion we had at the conference "Representations of p-adic Groups and the Langlands Correspondence, in honor of Colin Bushnell," held at King’s College London. It was through this conversation that I became aware of the significant work of Silberger and Zink, which greatly informed this research.
\newpage
\printnomenclature 
\newpage
\bibliographystyle{plain} 
\bibliography{references}

\end{document}